\documentstyle[amscd,amssymb,epic,eepic,12pt]{amsart}
\input xy
\xyoption{all}

\CompileMatrices                                                             
\newdir{((}{{\lhook\kern-.1em}}                                              
\newdir{))}{{\rhook\kern-.1em}}                                              
\newdir{ >}{{}*!/-5pt/@{>}}                                                  

\newcommand{\la}{\langle}
\newcommand{\ra}{\rangle}

\newtheorem{theorem}{\bf Theorem}[section]
\newtheorem{lemma}[theorem]{\bf Lemma}
\newtheorem{remark}[theorem]{\bf Remark}
\newtheorem{prop}[theorem]{\bf Proposition}
\newtheorem{corollary}[theorem]{\bf Corollary}
\newtheorem{definition}[theorem]{\bf Definition}

\newcommand{\be}{\begin{equation}}
\newcommand{\ee}{\end{equation}}        

\newfont{\bfc}{cmbsy10 scaled 1200}  
\newfont{\dr}{msbm10 scaled \magstep1}  
\newfont{\sdr}{msbm8}  
\newfont{\gl}{eufm10 scaled \magstep1}  

\DeclareFontFamily{OT1}{rsfs}{}
\DeclareFontShape{OT1}{rsfs}{n}{it}{<->rsfs10}{}
\DeclareMathAlphabet{\curly}{OT1}{rsfs}{n}{it}

\newcommand{\CC}{{\Bbb C}}

\newcommand{\PP}{{\Bbb P}}

\newcommand{\RR}{{\Bbb R}}
\newcommand{\ZZ}{{\Bbb Z}}

\newcommand{\glie}{{\frak g}}
\newcommand{\hlie}{{\frak h}}
\newcommand{\klie}{{\frak k}}
\newcommand{\llie}{{\frak l}}
\newcommand{\plie}{{\frak p}}

\newcommand{\tlie}{{\frak t}}
\newcommand{\ulie}{{\frak u}}
\newcommand{\zlie}{{\frak z}}

\newcommand{\fW}{{\frak W}}
\newcommand{\fV}{{\frak V}}
\newcommand{\fX}{{\curly X}}

\newcommand{\Ad}{\operatorname{Ad}}

\newcommand{\Chern}{{\frak C}}

\newcommand{\End}{\operatorname{End}}
\newcommand{\GL}{\operatorname{GL}}
\newcommand{\Gr}{\operatorname{Gr}}
\newcommand{\Hom}{\operatorname{Hom}}
\newcommand{\Id}{\operatorname{Id}}

\newcommand{\Ker}{\operatorname{Ker}}
\newcommand{\Lie}{\operatorname{Lie}}

\newcommand{\pt}{\operatorname{pt}}
\newcommand{\rk}{\operatorname{rk}}

\newcommand{\Tr}{\operatorname{Tr}}
\newcommand{\U}{\operatorname{U}}
\newcommand{\UU}{\operatorname{U(1)}}

\newcommand{\Vol}{\operatorname{Vol}}
\renewcommand{\exp}{\operatorname{exp}}

\newcommand{\ov}{\overline}

\newcommand{\plx}{{\pi^{\LLL}_X}}
\newcommand{\plf}{{\pi^{\LLL}_{\FFF}}}

\newcommand{\AAA}{{\curly A}}
\newcommand{\CCC}{{\curly C}}

\newcommand{\FFF}{{\cal F}}
\newcommand{\GGG}{{\curly G}}

\newcommand{\LLL}{{\curly L}}

\newcommand{\SSS}{{\curly S}}

\newcommand{\YMH}{{\cal YMH}}
\newcommand{\Met}{{\curly M}\!et^p_2}
\newcommand{\MetB}{{\curly M}\!et^p_{2,B}}
\newcommand{\imag}{i}
\newcommand{\qu}{/\kern-.7ex/}
\newcommand{\exh}{\to\kern-1.8ex\to}

\title[A Hitchin--Kobayashi correspondence for Kaehler fibrations]
{A Hitchin--Kobayashi correspondence for Kaehler fibrations}
\author{Ignasi Mundet i Riera}
\address{Departamento de Matem{\'a}ticas \\
Universidad Aut{\'o}noma de Madrid \\ 
Madrid, Spain}
\date{12 June 1999}
\email{ignasi.mundet@@uam.es}

\begin{document}
\maketitle
\begin{abstract}
Let $X$ be a compact Kaehler manifold and $E\to X$ a principal
$K$ bundle, where $K$ is a compact connected Lie group. 
Let $\AAA^{1,1}$ be the set of connections on $E$ whose
curvature lies in $\Omega^{1,1}(E\times_{\Ad}\klie)$. 
Let $\klie=\Lie(K)$, and fix on $\klie$ a nondegenerate biinvariant 
bilinear pairing. This allows to identify $\klie\simeq\klie^*$.
Let $F$ be a Kaehler left K-manifold and suppose that there exists
a moment map $\mu:F\to\klie^*$ for the action of $K$ on $F$.
Let $\SSS=\Gamma(E\times_K F)$. In this paper we study the equation 
$$\Lambda F_A+\mu(\Phi)=c$$
for $A\in\AAA^{1,1}$ on $E$ and a section $\Phi\in\SSS$,
where $F_A$ is the curvature of $A$ and $c\in\klie$ is a fixed central
element. We study which orbits
of the action of the complex gauge group on $\AAA^{1,1}\times\SSS$ 
contain solutions of the equation and we define a positive functional on 
$\AAA^{1,1}\times\SSS$ which generalises the Yang-Mills-Higgs functional
and whose local minima coincide with the solutions of the equation. 
\end{abstract}

\tableofcontents

\section{Introduction}
\label{intro}

\subsection{}
Let $X$ be a compact Kaehler manifold. 
Let $G$ be a connected complex reductive 
Lie group with maximal compact subgroup $K$, and let $E\to X$
be a $K$-principal bundle on $X$ (with the $K$ action on the right). 
Let $\GGG_K=\Gamma(E\times_{\Ad} K)$ be the real gauge group of $E$, and
let $\GGG_G=\Gamma(E\times_{\Ad} G)$ be the complex gauge group of $E$. 
$\GGG_G$ is the complexification of $\GGG_K$ and
is the gauge group of the $G$-principal bundle $E_G=E\times_K G$.
Let $\AAA$ be the space of $K$-connections on $E$.
The group $\GGG_K$ acts on $\AAA$ by pullback, and this action can
be extended to an action of $\GGG_G$ (see subsection \ref{setting}).
Let $\AAA^{1,1}\subset\AAA$ be the space of connections whose
curvature belongs to $\Omega^{1,1}(E\times_{\Ad}\klie)$
(equivalently, those which define an integrable holomorphic structure 
on $E_G$). The space $\AAA^{1,1}$ is $\GGG_G$-invariant. 

Let $F$ be any Kaehler manifold. Suppose that there is a Hamiltonian left 
action of $K$ on $F$ which respects the complex structure, and let 
$\mu:F\to\klie^*$ be a moment map for this action. We recall that by 
definition the following is satisfied:
(C1) for any $s\in\klie$, $d\mu(s)=\iota_{\fX_s}\omega_F$ (where 
$\fX_s$ is the field on $F$ generated by $s\in\klie$ and $\omega_F$
is the symplectic form of $F$) and
(C2) $\mu$ is equivariant with respect to the actions
of $K$ on $F$ and the coadjoint action on $\klie^*$.
The map $\mu$ is unique up to addition of constant central 
elements of $\klie^*$. 

Since $F$ is Kaehler, the action of $K$ on $F$ extends 
automatically to a unique holomorphic action of 
$G$ (see \cite{GS}). Let $\FFF=E\times_K F=E_G\times_G F\to X$ be 
the associated bundle on $X$ with fibre $F$, and let
$\SSS$ be the space $\Gamma(\FFF)$ of smooth sections of $\FFF$.
The group $\GGG_G$ acts on $\FFF$, and consequently also on $\SSS$.
Since $\mu$ is $K$-equivariant 
we can extend fibrewise the moment map $\mu$, thus obtaining
for any $\Phi\in\SSS$ a section $\mu(\Phi)\in \Omega^0(E\times_{\Ad}\klie^*)$.

In this paper we study the equation
\begin{equation}
\Lambda F_A+\mu(\Phi)=c,
\label{equacio}
\end{equation}
where $A\in\AAA^{1,1}$, $F_A\in\Omega^2(E\times_{\Ad}\klie)$ is 
the curvature of $A$,
$\Phi\in\SSS$ and $c\in\Omega^0(E\times_{\Ad}\klie)$
is a constant central element. Here $\Lambda:\Omega^*(X)\to\Omega^{*-2}(X)$ 
is the adjoint of the map given by wedging with the symplectic form 
$\omega$ of $X$, and we identify (by means of a biinvariant metric on 
$\klie$) $\Omega^0(E\times_{\Ad}\klie^*)$ with 
$\Omega^0(E\times_{\Ad}\klie)$.

\subsection{}
The main question which we consider is the following: for which pairs
$(A,\Phi)\in \AAA^{1,1}\times\SSS$ there exist a gauge
transformation $g\in\GGG_G$ such that
$(A',\Phi')=g(A,\Phi)$ satisfies equation (\ref{equacio})?
We will define two conditions on pairs $(A,\Phi)$ called simplicity and
$c$-stability, and in Theorem \ref{main} 
we will prove that, if $(A,\Phi)$ is a simple pair,
then there exist a gauge $g\in\GGG_G$ sending
$(A,\Phi)$ to a pair $g(A,\Phi)$ which solves (\ref{equacio})
if and only if $(A,\Phi)$ is $c$-stable.
Observe that if $g(A,\Phi)$ solves (\ref{equacio}), so does
$kg(A,\Phi)$ for any $k\in\GGG_K$. We will also prove that
in each $\GGG_G$ orbit inside $\AAA^{1,1}\times\SSS$ there
is at most one $\GGG_K$ orbit of pairs which satisfy (\ref{equacio}).
This is proved in Theorem \ref{main}.
Such a characterization of solutions to (\ref{equacio})
is typically called a Hitchin--Kobayashi correspondence,
since a particular case of it ($F$ equal to a point) was independently
conjectured by Hitchin and Kobayashi.

One can look at Theorem \ref{main} from two different points of view.
When $X$ consists of a single point, the curvature term vanishes
in equation (\ref{equacio}), and so our problem reduces to a well
known one in Kaehler geometry. Namely, that of studying which $G$ orbits
inside $F$ contain zeroes of the moment map $\mu$.
More generally, one studies which $G$ orbits have points whose image is
a fixed central element in $\klie^*$ or belongs to a given coadjoint orbit
in $\klie^*$. If $F$ is a projective manifold,
one can answer as follows:
a $G$ orbit contains a zero of the moment map if and only if
it is stable in the sense of Mumford Geometric Invariant Theory 
(GIT for short) \cite{KeNe, MFK, GS}. 
To extend the notion of GIT stability to actions on any Kaehler manifold
$F$, we use the notion of analytic stability (see definition
\ref{estabilitat}).
This notion coincides with that of GIT stability in the case of projective 
manifolds, and characterizes the $G$-orbits 
in which the moment map vanishes somewhere (see Theorem \ref{corr}).
This is the content of the so-called Kempf--Ness theory. 
So, in this sense, our result can be viewed as a fibrewise
generalisation of Kempf--Ness theory. 

There is, however, another point of view which allows to look at 
Theorem \ref{main} as a result {\it {\`a} la} Kempf--Ness in infinite 
dimensions. One can give a Kaehler structure to the configuration space
$\AAA^{1,1}\times\SSS$ (for this we use the same biinvariant metric on
$\klie$ that was used to give a sense to equation (\ref{equacio})); 
then the action of the gauge group
$\GGG_K$ is symplectic and by isometries, and the left hand side in
equation (\ref{equacio}) is a moment map of this action 
(see sections \ref{exemp1}, \ref{exemp2} and \ref{disgressio}). This point 
of view was adopted for the first time in the context of gauge theories
by Atiyah and Bott \cite{AB} in their study of Yang-Mills equations
over Riemann surfaces, which are a particular case of the equations
that we consider. The idea of Atiyah and Bott was used by Donaldson
\cite{Do1} in his proof of the theorem of Narasimhan and Seshadri
(which is a particular case of Theorem \ref{main}), and it has been 
subsequently often used in studying other particular cases of equation 
(\ref{equacio}).

\subsection{}
After proving Theorem \ref{main} we address the problem of finding
a functional on $\AAA^{1,1}\times\SSS$ which generalises the
classical Yang-Mills-Higgs functional and whose (local) minima
satisfy equation (\ref{equacio}). 
We define for any connection $A$ on $E$ a {\it covariant derivation}
which assigns to any section $\Phi\in\SSS$ a section
$d_A\Phi\in\Omega^1(\Phi^*\Ker d\pi_F)$, where 
$\pi_F:\FFF\to X$ denotes the projection.
When $F$ is a vector space on which $K$ acts linearly, $\FFF$
is a vector bundle, $\Ker d\pi_F$ is canonically isomorphic to $\FFF$, 
and the covariant derivation $d_A$ coincides with
the usual one in differential geometry.
The Yang-Mills-Higgs functional is defined as 
$$\YMH_c(A,\Phi)=\|F_A\|_{L^2}^2+\|d_A\Phi\|_{L^2}^2
+\|c-\mu(\Phi)\|_{L^2}^2,$$
where $(A,\Phi)\in\AAA^{1,1}\times\SSS$. 
If $F$ is a representation space for $K$, the Yang-Mills-Higgs
functional coincides with the usual one in gauge theories.
Now, using the splitting 
$\Omega^1(X)\otimes\CC=\Omega^{1,0}(X)\oplus\Omega^{0,1}(X)$
we obtain from $d_A$ an operator $\overline{\partial}_A$
which sends any $\Phi\in\SSS$ to a section 
$\overline{\partial}_A\Phi\in\Omega^{0,1}(\Phi^*\Ker d\pi_F)$.
We then consider the two equations for a connection $A\in\AAA^{1,1}$ 
and a section $\Phi\in\SSS$
\begin{equation}
\left\{\begin{array}{l}
\overline{\partial}_A\Phi=0, \\
\Lambda F_A+\mu(\Phi)=c.\end{array}\right.
\label{equacions}
\end{equation}
We show in section \ref{YMHs} that the pairs 
$(A,\Phi)\in\AAA^{1,1}\times\SSS$ solving these 
equations minimize the Yang-Mill-Higgs functional among the pairs whose
section belong to a fixed homology class of sections of $\FFF$.

The way we identify the solutions of the equations with 
(local) minima of Yang-Mills-Higgs functional is similar to the 
one used in the study of holomorphic pairs (see \cite{Br1}). The main
difference is in the step where in dealing with holomorphic
pairs \cite{Br1} uses the Kaehler identities. At that point we use 
certain results on the coupling form on symplectic fibrations
due to Guillemin, Lerman and Sternbert \cite{GLeS}.

Note that in the Hitchin--Kobayashi correspondence we ignore the
first equation in (\ref{equacions}), that is, $\overline{\partial}_A\Phi=0$
(this equation can be given a sense even when $F$ is not a vector
space; see section \ref{YMHs}). Indeed, this condition is not necessary
in the proof: we only need $\Phi$ to be smooth. Furthermore,
the equation $\overline{\partial}_A\Phi=0$ is invariant under
$\GGG_G$, while the interest of our problem stems from the fact
that $\Lambda F_A+\mu(\Phi)=c$ is only $\GGG_K$-invariant.

It is a remarkable fact that both the equations (\ref{equacions}) and the 
results in section \ref{YMHs} make perfect sense even when the complex
structure on $F$ is not integrable. In a forthcoming paper we will study these
equations and we will show how the gauge equivalence classes of its solutions
can be used to define invariants of Hamiltonian actions on compact
symplectic manifold (see \cite{Mu}).

\subsection{}
Many particular instances of equations (\ref{equacions}) have been 
already studied. When $F=\{\pt\}$ equation (\ref{equacio}) becomes the 
Hermite--Einstein equation, which was studied for example in
\cite{BarTi, Do1, Do2, NSe, UY}. A good reference for Hitchin--Kobayashi
correspondence for Hermite--Einstein equations and its interesting
history is the book by L{\"u}bke and Teleman \cite{LTe}.
When $F$ is a representation space for $K$, the fibre bundle $\FFF$ is a 
vector bundle. Theorem \ref{main} has been proved for many particular choices 
of $K$ and representations $K\to\U(F)$ (see for example 
\cite{Br1, Br2, BrGP3, GP1, GP2, Hi, JT, Si}). In 1996 Banfield \cite{Ba} 
gave a proof of the Hitchin--Kobayashi correspondence for any $K$ and any 
representation space $F$ of $K$.

A particular case of our construction which does not fit
in Banfield's result is that of extensions and filtrations of 
vector bundles.
They arise when $F$ is a Grassmannian or, more generaly,
any flag manifold. A Hitchin--Kobayashi correspondence for
extensions was studied by Bradlow and
Garc{\'\i}a--Prada \cite{BrGP1}, and by Daskalopoulos, Uhlenbeck 
and Wentworth \cite{DaUW}; the correspondence for filtrations 
has been proved by {\'A}lvarez C{\'o}nsul and Garc{\'\i}a--Prada 
\cite{AlGP}.

\subsection{}
This paper is organised as follows. 
In section \ref{statement} we state the main
result of this paper. Sections \ref{integral} to \ref{demostracio}
are devoted to the proof of this result. In section \ref{integral}
we explain the construction of a certain functional which will be
the main tool in the proof. In section \ref{esKaehler} we describe
a Kaehler structure on the manifold $\AAA^{1,1}\times\SSS$ and
we identify our equation as a moment map for the action of
$\GGG_K$ on $\AAA^{1,1}\times\SSS$. In section \ref{correspondencia}
we prove a particular case of our theorem, and the general proof
is given in section \ref{demostracio}. In section \ref{YMHs}
we introduce (a generalisation of) the Yang-Mills-Higgs
functional and we prove that its minima coincide with the solutions
of equations (\ref{equacions}). Finally, in sections
\ref{banfield} and \ref{filtracions} we work out two different
examples of our correspondence.

\subsection{Acknowledgements} This paper is part of my Ph. D. Thesis.
It is for me a pleasure to thank my advisor, Oscar Garc{\'\i}a--Prada,
for his continuous support and encouragement, and for his excellent
guidance. I also thank Vicente Mu{\~n}oz for carefully reading 
this paper and for his useful comments. Finally, I thank the referee
for his clarifying observations and especially for pointing out the reference
\cite{GLeS}, which has greatly simplified section \ref{YMHs}.

\section{Stability and statement of the correspondence}
\label{statement}

\subsection{The isomorphism $\klie\simeq\klie^*$}
\label{isotriat}
To give a meaning to equation (\ref{equacio}) we need a $K$-equivariant
isomorphism $\klie\simeq\klie^*$. From now on we will assume that such 
an isomorphism comes from a biinvariant metric on $K$ which is the 
pullback of the Killing metric through a faithful representation
$\rho_a:K\to U(W_a)$, where $W_a$ is a Hermitian vector space. 
(In other words, the isomorphism $\klie\simeq\klie^*$ is a {\it hidden} 
parameter of the equation, and we prove the correspondence for some 
particular choices of it.)
Our characterisation of solutions to (\ref{equacio}) will depend on the
choice of $\rho_a$ and $W_a$ (this is not strange, since the equation
also depends on them). 

\subsection{The action of $\GGG_G$ on $\AAA$}
\label{setting}
Let $\omega$ be the symplectic form on $X$ and $I\in\End(TX)$ 
the complex structure. In the sequel $\omega^{[k]}$ will denote 
$\omega^k/k!$ The volume element $\omega^{[n]}$ will be implicitly 
assumed in all the integrals of functions on $X$.

Let $\CCC$ be the set of $G$-invariant complex
structures on $E_G=E\times_K G$ for which the map
$d\pi_G:TE_G\to \pi_G^*TX$ is complex. 
We define a map $\Chern:\CCC\to\AAA$
as follows. A complex structure $I\in\CCC$ is mapped to the
connection $\Chern(I)$ given by the horizontal distribuition
$I(TE)\cap TE\subset TE$ (this makes sense, since the inclusion
$E=E\times_K K\subset E\times_K G$ given by $K\subset G$ induces an inclusion
$TE\subset TE_G$). This defines a connection and the map $\Chern$ is a 
bijection (see \cite{Sn}). We call $\Chern$ the Chern map.

The following is readily checked.

\begin{lemma}
(i) Let $G$ act holomorphically on a vector space $W$. Let us take
a complex structure $I\in\CCC$.
The associated bundle $V=E_G\times_G W$ is endowed by $I$ of a complex
structure $I_V\in\End(TV)$.
Any section $\sigma\in\Omega^0(V)$ may be viewed as a map
$\sigma:X\to V$. Then, the antiholomorphic part 
$\overline{\partial}_I(\sigma)
=(d\sigma+I_V\circ d\sigma\circ I_X)/2$ can be  
regarded as an element in $\Omega^{0,1}(V)$.

(ii) For any $A\in\AAA$ we have
$\overline{\partial}_A=\overline{\partial}_{\Chern^{-1}(A)}$,
where $\overline{\partial}_A:\Omega^0(V)\to\Omega^{0,1}(V)$
is the usual $\overline{\partial}$ operator obtained from $A$.

(iii) The set $\AAA^{1,1}$ is mapped by $\Chern^{-1}$ to the set of 
integrable complex structures on $E_G$. 
\label{integr}
\end{lemma}

The group $\GGG_G$ acts on $\CCC$ by pullback, and 
using the map $\Chern$ we transfer the action of $\GGG_G$ on $\CCC$
to an action on $\AAA$. This action extends the action of $\GGG_K$ and
(by {\it (iii)} in the preceeding lemma) leaves invariant the 
subset $\AAA^{1,1}\subset\AAA$.

\subsection{Maximal weights}
\label{grkaehler}
Let $I_F\in\End(TF)$ be the complex structure of $F$. 
We will denote by $\la u,v\ra=\omega_F(u,I_Fv)$ the Kaehler metric on $F$.

Let $s\in\klie$ be any nonzero element, and let us write 
$\mu_{s}=\la\mu,s\ra_{\klie}:F\to\RR$.
(Here and in the sequel we denote by $\la\cdot,\cdot\ra_W:W^*\times W\to\RR$ 
the canonical pairing for any vector space $W$.)
Recall that $\fX_s$ is the field generated on $F$ by $s$.

\begin{lemma}
The gradient of $\mu_{s}$ is $I_F\fX_{s}$.
\label{gradient}
\end{lemma}
\begin{pf}
Let $x\in F$ and take any vector $v\in T_xF$.
Then $\nabla_v(\mu_{s})=\la d\mu_{s},v\ra_{T_xF}=
\omega_F(\fX_{s},v)=\omega_F(I_F\fX_{s},I_Fv)=\la I_F\fX_{s},v\ra,$ 
by the definition of moment map. 
\end{pf}

Consider the gradient flow $\phi^t_{s}:F\to F$ of the function
$\mu_{s}$. $\phi^t_s$ is defined by these properties: 
$\phi^0_{s}=\Id$ and $\frac{\partial}{\partial t}\phi^t_s
=\nabla(\mu_{s})=I\fX_{s}$.
Using the action of $G$ on $F$ we can write $\phi^t_{s}(x)=e^{\imag ts}x$.

\begin{definition}
Let $x\in F$ be any point, and take an element $s\in\klie$.
Let $$\lambda_t(x;s)=\mu_s(e^{\imag ts}x).$$
We define the {\rm maximal weight $\lambda(x;s)$ of the action 
of $s$ on $x$} to be
$$\lambda(x;s)=\lim_{t\to \infty}\lambda_t(x;s)
\in\RR\cup\{\infty\}.$$
\label{defpesmax}
\end{definition}

This limit always exists since by Lemma \ref{gradient} the function
$\lambda_t(x;s)$ increases with $t$.
The definition of the maximal weight depends on the
chosen moment map. Since this is not unique, we will sometimes 
write the maximal weight of $s\in\klie$ acting on $x\in F$
with respect to the moment map $\mu$ as $\lambda^{\mu}(x;s)$.

\begin{prop}
The maximal weights satisfy the following properties:
\begin{enumerate}
\item They are $K$-equivariant, that is, for any $k\in K$,
$\lambda(kx;ksk^{-1})=\lambda(x;s)$.
\item For any positive real number $t$ one has
$\lambda(x;ts)=t\lambda(x;s)$.
\end{enumerate}
\label{proppes}
\end{prop}

See sections \ref{banfield} and \ref{filtracions} for explicit
computations of maximal weights in some particular situations.

\subsection{Parabolic subgroups}
\label{parabolic}
A good reference for this material is \cite{R}.
Let $\glie$ be the Lie algebra of $G$, and split $\glie=\zlie\oplus
\glie^s$ as the sum of the centre plus the semisimple part 
$\glie^s=[\glie,\glie]$ of $\glie$.
Take a Cartan subalgebra $\hlie\subset\glie^s$. 
Let $R\subset\hlie^*$ be the set of roots.
We can decompose
$$\glie=\zlie\oplus\hlie\oplus\bigoplus_{\alpha\in R}\glie_{\alpha},$$
where $\glie_{\alpha}\subset\glie^s$ is the subspace on which
$\hlie$ acts through the character $\alpha\in\hlie^*$.

Fixing a (irrational) linear form on $\hlie^*$, we divide the set
of roots in positive and negative roots: $R=R^+\cup R^-$.
Let us write the set of simple roots $\Delta=(\alpha_1,\dots,\alpha_r)
\subset R^+$. Recall that the set $\Delta$ is characterised by the
following property: any root can be written as 
a linear combination of the elements of $\Delta$ with integer
coefficients all of the same sign. Furthermore, $r$ equals 
$\dim_{\CC}\hlie$, the rank of $G$. The simple coroots are by definition 
$\alpha_j'=2\alpha_j/\la\alpha_j,\alpha_j\ra$, where $1\leq j\leq r$.

We have taken a maximal compact subgroup $K\subset G$.
From now on we will assume that the following relation holds between $K$ and
the Cartan subalgebra $\hlie$: $\zlie\oplus\hlie$ is the
complexification of the Lie algebra $\tlie$ of a maximal torus $T\subset K$.

\begin{lemma}
Chose, for any root $\alpha\in R$, a nonzero element
$g_{\alpha}\in\glie_{\alpha}$ in such a way that $g_{\alpha}$
and $g_{-\alpha}$ satisfy $\la g_{\alpha},g_{-\alpha}\ra=1$.
Let $\RR R^*\subset \hlie$ denote the real span of the 
duals (with respect to the Killing metric) of the roots. 
Assume that $\zlie\oplus\hlie$ is the complexification of the Lie algebra
of a maximal torus $T$ of a maximal compact subgroup $K\subset G$.
Then $\glie^s\cap\klie=\imag\RR R^*\oplus \bigoplus_{\pm\alpha\in R}
\RR(g_{\alpha}+g_{-\alpha})\oplus \RR(\imag g_{\alpha}-\imag g_{-\alpha}).$
\label{compact}
\end{lemma}

This lemma (and the following ones in this subsection) can be easily proved 
using basic results on reductive Lie groups (see for example \cite{FH}).

Let $\lambda_1,\dots,\lambda_r$ be the set of fundamental weights,
which belong to $\hlie^*$ and are the duals with respect to the 
Killing metric of the simple coroots.
Let us denote by $\lambda'_1,\dots,\lambda'_r$ the elements in $\hlie$
dual to the fundamental weights through the Killing metric.

To define a parabolic subgroup of $G$, take any subset
$A=\{\alpha_{i_1},\dots,\alpha_{i_s}\}\subset \Delta$. Let
$$D=D_A=\{\alpha\in R \mid \alpha=\sum_{j=1}^r m_j\alpha_j
\mbox{, where $m_{i_t}\geq 0$ for $1\leq t\leq s$}\}.$$

\begin{definition}
The subalgebra $\plie=\zlie\oplus
\hlie\oplus\bigoplus_{\alpha\in D}\glie_{\alpha}$
will be called the {\rm parabolic subalgebra} of $\glie$ 
with respect to the set $A\subset\Delta$. 
The connected subgroup $P$ of $G$ whose subalgebra is $\plie$ will be called 
the {\rm parabolic subgroup} of $G$ with respect to $A$.
Furthermore, any positive (resp. negative) linear combination of the 
fundamental weights $\lambda_{i_1},\dots,\lambda_{i_s}$ plus an element
of the dual of $\imag(\zlie\cap\klie)$ will be called 
a {\rm dominant} (resp. {\rm antidominant}) {\rm character} on 
$\plie$ (or on $P$).
\end{definition}

\begin{remark}
We will regard $G$ as a parabolic subgroup 
of itself (with respect to the empty set $\emptyset\subset \Delta$).
\end{remark}

Observe that our definition of parabolic subgroup depends upon the 
choice of a Cartan subalgebra $\hlie\subset\glie$ and of a linear
form on $\hlie^*$. In general, any parabolic subgroup
$P\subset G$ obtained from a different choice of Cartan subalgebra
and linear form will be conjugate to a parabolic subgroup obtained
from our data.

\subsection{Parabolic subgroups and filtrations}
\label{parfil}
Let $\rho:K\to U(W_{\rho})$ be a representation on a Hermitian vector space
$W_{\rho}$.
We will write its (unique) lift to a holomorphic representation
of the complexification $G$ of $K$
with the same letter $\rho:G\to GL(W_{\rho})$. Take $P\subset G$ to be 
the parabolic subgroup with respect to a set $A=\{\alpha_{i_1},\dots,
\alpha_{i_s}\}\subset\Delta$. Let $\chi$ be the dual of 
an antidominant character of $P$.
Thanks to our conventions (Lemma \ref{compact}), $\chi$ belongs
to $\imag\klie$. So, since $\rho$ is unitary, $\rho(\chi)$ diagonalises
and has real eigenvalues.
Let $\lambda_1<\dots<\lambda_r$ be the set of different
eigenvalues of $\rho(\chi)$, and let us write $W(\lambda)$ the eigenspace
of eigenvalue $\lambda$. 
Let $W^{\lambda_k}=\bigoplus_{j\leq k}W(\lambda_j)$, 
and let $\fW_{\rho}(\chi)$ be the partial flag
$0\subset W^{\lambda_1}\subset\dots\subset W^{\lambda_r}=W_{\rho}$.

\begin{lemma}
(i) The action of $P$ leaves invariant the partial flag $\fW_{\rho}(\chi)$.
Suppose that the restriction of $\rho$ to the semisimple part
$\plie^s$ of $\plie$ is faithful.
If $\chi=z+\sum_{k=1}^s m_k\lambda'_{i_k}$, where $z\in\zlie$, and, for any
$k$, $m_k<0$, then $P$ is precisely the antiimage by $\rho$
of the stabiliser of $\fW_{\rho}(\chi)$.
(ii) Let $\chi\in\imag\klie$ be any element. There is a choice of 
Cartan subalgebra $\hlie\subset\glie$ contained in $\plie$
such that $\chi\in\hlie$ and $\chi$ is antidominant with 
respect to $P$ if and only if the stabiliser of the partial
flag $\fW_{\rho}(\chi)$ contains $P$.
\label{equi1}
\end{lemma}

\begin{lemma}
Let $\chi$ be any element in $\imag\klie$. The antiimage by 
$\rho$ of the stabiliser of $\fW_{\rho}(\chi)$ is a parabolic
subgroup $P_{\rho}(\chi)$ of $G$. Moreover, $\chi$ 
is the dual of an antidominant character of $P_{\rho}(\chi)$.
\label{equi3}
\end{lemma}

Let us take now any subspace $W'\subset W_{\rho}$ 
belonging to the filtration $\fW_{\rho}(\chi)$. 
Since $P$ leaves $W'$ invariant, we may define 
$\overline{W}'= G\times_P W'\to G/P$ (here we view $G$ as a right 
$P$-principal bundle). Define also an action of
$G$ on $G\times W'$ by $g'(g,w)=(g'g,g^{-1}g'gw)$. This action 
descends to an action on $\overline{W}'$.
Repeating this for each subspace in $\fW_{\rho}(\chi)$ we obtaing the
following.

\begin{lemma}
The filtration of holomorphic vector bundles 
$\overline{\fW}_{\rho}(\chi)=G\times_{\rho}\fW_{\rho}(\chi)\to G/P$
admits a holomorphic lift of the right action of $G$ on $G/P$.
\label{tautologic}
\end{lemma}

\subsection{Parabolic and maximal compact subgroups}
\label{maxcomp}
Given any parabolic subgroup $P\subset G$ with Lie algebra $\plie$,
we will write $P_K$ (resp. $\plie_K$) for the subgroup
$P\cap K$ (resp. the subalgebra $\plie\cap\klie$). 
$P_K$ is a maximal compact subgroup of $P$.

\begin{lemma}
Let $E_G\to X$ be a $G$-principal bundle on any topological space $X$.
If $E_G$ admits reductions of its structure group from 
$G$ to a parabolic subgroup $P$ and to the maximal compact subgroup $K$,
then it also admits a reduction of its structure group from $G$
to $P_K$.
\label{maxcomp1}
\end{lemma}

\begin{lemma}
Let $P$ be a parabolic subgroup with respect to the set
$$A=\{\alpha_{i_1},\dots,\alpha_{i_s}\}\subset\Delta.$$
For any $j\in\{i_1,\dots,i_s\}$, the element $\lambda'_j\in\imag\klie$ 
(dual with respect to the Killing metric of the fundamental weight
$\lambda_j$) is left fixed by the adjoint action of $\plie_K$
on $\glie$.
\label{maxcomp2}
\end{lemma}

\subsection{Reductions of the structure group and filtrations}
\label{reduccions}
Following the notation in subsection \ref{parfil}, we denote
$V_{\rho}=E\times_{\rho}W_{\rho}$. In this subsection we will see that 
there is a correspondence between the 
reductions of the structure group of $E$ to a parabolic subgroup $P$ 
together with an antidominant character of $P$, and certain filtrations of 
$V_{\rho}$ by subbundles. We denote $E(G/P)$ the bundle $E_G\times_G(G/P)$.
The space of reductions of the structure group of $E_G$ from $G$
to $P$ is $\Gamma(E(G/P))$. 

\subsubsection{}
\label{sectfilt}
Fix a parabolic subgroup $P\subset G$ and take a reduction
$\sigma\in\Gamma(E(G/P))$. Let $\chi$ be an antidominant character
for $P$. There is a canonical reduction of
the structure group $G$ of $E_G$ to $K$, since $E_G=E\times_K G$.
Thanks to Lemma \ref{maxcomp1}, this reduction, together
with $\sigma$, gives a reduction $\sigma_K\in \Gamma(E(G/P_K))$, where
$P_K=P\cap K$. And then, Lemma \ref{maxcomp2} implies that
we get a section $g_{\sigma,\chi}\in\Omega^0(E\times_{\Ad}\imag\klie)
=\imag\Lie(\GGG_K)$ which is fibrewise the dual of $\chi$.

With the element $g_{\sigma,\chi}$ we can obtain a filtration of
$V_{\rho}$ as follows. First of all, $\rho(g_{\sigma,\chi})$
has constant real eigenvalues (which are equal to those of
$\rho(\chi)\in\End(W_{\rho})$). Let $\lambda_1<\dots<\lambda_r$ 
be the different eigenvalues, and let $V_{\rho}(\lambda_j)$ 
be the eigenbundle of eigenvalue $\lambda_j$.
Finally, let $V_{\rho}^{\lambda_k}=\bigoplus_{i\leq k}V_{\rho}(\lambda_j)$.
Denote by $\fV_{\rho}(\sigma,\chi)$ the filtration 
$$0\subset V_{\rho}^{\lambda_1}\subset
V_{\rho}^{\lambda_2}\subset\dots\subset V_{\rho}^{\lambda_r}=V_{\rho}.$$

Alternatively, recall that on $G/P$ there is a 
filtration of $G$-equivariant (holomorphic)
vector bundles, $\overline{\fW}_{\rho}(\chi)$ (see Lemma \ref{tautologic}). 
$G$-equivariance allows to define the filtration
$\overline{\fV}_{\rho}(\chi)=E\times_G \overline{\fW}_{\rho}(\chi)\to E(G/P)$.
Then $\fV_{\rho}(\sigma,\chi)=\sigma^*\overline{\fV}_{\rho}(\chi)$.

\subsubsection{}
\label{filtsect}
Conversely, take $g\in\Omega^0(E\times_{\Ad}\imag\klie)$.
Suppose that $\rho(g)$ has constant eigenvalues, and let
$\lambda_1<\dots<\lambda_r$ be the set of different values
they take. Just as before, we consider the filtration
\begin{equation}
0\subset V_{\rho}^{\lambda_1}\subset
V_{\rho}^{\lambda_2}\subset\dots\subset V_{\rho}^{\lambda_r}=V_{\rho}.
\label{asterix}
\end{equation}
Fix a point $x\in X$. After trivialising the fibre $E_x$ we
can identify $g(x)$ with and element $\chi$ of $\imag\klie$.
Let $P=P_{\rho}(\chi)$ (see Lemma \ref{equi3}).
We obtain a reduction $\sigma\in\Gamma(E(G/P))$ as follows.
Let $y\in X$. Trivialise $E_y$ and identify $g(y)$ with
$\chi_y\in\imag\klie$. Let 
$$\sigma(y)=\{g\in G\mid g(\fW_{\rho}(\chi))=\fW_{\rho}(\chi_y)\}.$$
Then $\sigma(y)$ is invariant under left multiplication by
elements of $P$, and in fact gives a unique point in $G/P$
(here we use Lemma \ref{equi3}). Furthermore, the definition of 
$\sigma(y)$ is compatible with change of trivialisation in the
sense that it gives a section $\sigma\in\Gamma(E(G/P))$.

\begin{lemma}
The filtration (\ref{asterix}) is equal to $\fV_{\rho}(\sigma,\chi)$.
\end{lemma}

\subsubsection{Holomorphic reductions of the structure group}
\label{holomred}
Suppose that there is a fixed (integrable) holomorphic structure on $E_G$.
This structure induces a holomorphic structure
on the total space of the associated bundle $E(G/P)$, since
$G/P$ is a complex manifold and the action of $G$ on $G/P$
is holomorphic. 

\begin{definition}
Let $\sigma\in \Gamma(E(G/P))$. A reduction $\sigma$
is {\rm holomorphic} if the map $\sigma:X\to E(G/P)$
is holomorphic.
\end{definition}

One can give an equivalent definition of holomorphicity in terms
of the filtrations induced by the reduction $\sigma$ in the 
associated vector bundles.

\begin{lemma}
Let $\sigma\in \Gamma(E(G/P))$. If the reduction $\sigma$ is holomorphic
then, for any antidominant character $\chi$ of $P$ and for any representation
$\rho:K\to U(W)$, the filtration $\fV_{\rho}(\sigma,\chi)$ of $V_{\rho}$ is
holomorphic. Conversely, let $g\in\Omega^0(E\times_{\Ad}\imag\klie)$ have 
constant eigenvalues, and let $P\subset G$, $\sigma\in\Gamma(E(G/P))$, 
$\chi\in\imag\klie$ and $\fV_{\rho}(\sigma,\chi)$ be obtained from it 
as in \ref{filtsect}. Suppose that $\rho$ is faithful.
If $\fV_{\rho}(\sigma,\chi)$ is holomorphic, then so is $\sigma$.
\end{lemma}

\subsection{Total degree of a reduction of the structure group}
Let $V=V_{\rho_a}=E\times_{\rho_a}W_a$ be the vector bundle
associated to the representation $\rho_a$ (see subsection \ref{isotriat}).
We will apply the preceeding results on filtrations of vector bundles to $V$.
Let $P$ be a parabolic subgroup of $G$ with respect to 
$\{\alpha_{i_1},\dots,\alpha_{i_s}\}\subset \Delta$. Suppose that 
$\sigma\in \Gamma(E(G/P))$ is a reduction. Let $\chi$ be an antidominant
character of $P$.

We begin by defining the degree of the pair $(\sigma,\chi)$.
Let $0\subset V^{\lambda_1}\subset\dots\subset V^{\lambda_r}=V$ be
the filtration $\fV_{\rho_a}(\sigma,\chi)$ of $V$. 
For any vector bundle $V'$ we denote
$$\deg(V')={2\pi}\la c_1(V')\cup [\omega^{[n-1]}],[X]\ra.$$
Here $[\omega^{[n-1]}]$ denotes the cohomology class 
represented by the form $\omega^{[n-1]}$ and $[X]\in H_{2n}(X;\ZZ)$
is the fundamental class of $X$. Then we set
$$\deg(\sigma,\chi)=\lambda_r\deg(V)+
\sum_{k=1}^{r-1}(\lambda_k-\lambda_{k+1})\deg(V^{\lambda_k}).$$

\subsection{Stability, simple pairs and the correspondence}
\label{ssp}
Let $\sigma\in \Gamma(E(G/P))$ be a reduction.
We define the {\rm maximal weight} of $(\sigma,\chi)$ acting on a section
$\Phi\in\SSS$ as
$$\int_{x\in X}\lambda(\Phi(x);-\imag g_{\sigma,\chi}(x)),$$
where $\lambda(\Phi(x);-g_{\sigma,\chi}(x))$ is the maximal
weight of $-g_{\sigma,\chi}(x)$ acting on $\Phi(x)$ as defined
in \ref{defpesmax} (note that here we use the $K$-equivariance of
the maximal weights, as stated in Lemma \ref{proppes}).

Finally, given any central element $c\in \zlie\cap\klie$
we define the {\rm $c$-total degree} of the pair $(\sigma,\chi)$
as $$T^c_{\Phi}(\sigma,\chi)=\deg(\sigma,\chi)
+\int_{x\in X}\lambda(\Phi(x);-\imag g_{\sigma,\chi}(x))
+\la\imag\chi,c\ra\Vol(X).$$
Just as the maximal weights, the $c$-total degree is allowed
to be equal to $\infty$.

Now suppose that $X_0\subset X$ has as complement in $X$ a complex 
codimension 2 submanifold. Suppose also that a reduction $\sigma$
is defined only in $X_0$, that is, $\sigma\in\Gamma(X_0;E(G/P))$.
In this case it also makes sense to speak about 
$T^c_{\Phi}(\sigma,\chi)$ for any antidominant character $\chi$.
The only difficulty would be in defining the degree $\deg(\sigma,\chi)$;
however, it is well known that the degree of a vector bundle
can be computed by integrating the Chern-Weil form in the complement 
of a complex codimension 2 variety.

\begin{definition}
A pair $(A,\Phi)\in\AAA^{1,1}\times\SSS$ is {\rm $c$-stable} if
for any $X_0\subset X$ whose complement on $X$ is a complex 
codimension 2 submanifold, for any parabolic subgroup $P$ of $G$, 
for any holomorphic (with respect to the complex structure
$\Chern^{-1}A$ on $E_G$, see Lemma \ref{integr}) reduction
$\sigma\in \Gamma(X_0;E(G/P))$ defined on $X_0$, and for any antidominant 
character $\chi$ of $P$ we have $$T^c_{\Phi}(\sigma,\chi)>0.$$
\label{parella_estable}
\end{definition}

We will say that an element $s\in\GGG_G$ is {\rm semisimple} if,
for any $x\in X$, after identifying $(E\times_{\Ad}\glie)_x\simeq \glie$,
$s(x)\in \glie$ is a semisimple element.
(This is independent of the chosen isomorphism  
$(E\times_{\Ad}\glie)_x\simeq \glie$,
because an element of $\glie$ is semisimple if and only if any element 
in its orbit by the adjoint action of $G$ on $\glie$ is semisimple.)

\begin{definition}
A pair $(A,\Phi)$ is {\rm simple} if no semisimple element in 
$\Lie(\GGG_G)$ leaves $(A,\Phi)$ fixed,
that is, for any semisimple $s\in\Lie(\GGG_G)$, 
$\fX^{\AAA\times\SSS}_s(A,\Phi)\neq 0$.
\label{defsimple}
\end{definition}

\begin{remark}
If $(A,\Phi)$ is simple then so is any point in
the $\GGG_G$ orbit through $(A,\Phi)$.
\end{remark}

We are now ready to state the main theorem of this paper.

\begin{theorem}[Hitchin--Kobayashi correspondence]
Let $(A,\Phi)\in\AAA^{1,1}\times\SSS$ be a simple pair.
There exists a gauge transformation $g\in\GGG_G$ such that
\begin{equation}
\Lambda F_{g(A)}+\mu(g(\Phi))=c 
\label{hk}
\end{equation}
if and only if $(A,\Phi)$ is $c$-stable. Furthermore,
if two different $g,g'\in\GGG_G$ solve equation (\ref{hk}), then
there exists $k\in\GGG_K$ such that $g'=kg$.
\label{main}
\end{theorem}

We briefly explain the idea of the proof of Theorem \ref{main}.
We construct on $\AAA^{1,1}\times\SSS\times\GGG_G$ a functional $\Psi$
(that we will call integral of the moment map) whose critical
points give the solutions of equation (\ref{hk}). 
We prove that the pair $(A,\Phi)$ is $c$-stable if and only if the 
functional $\Psi$ is, in a certain sense, proper along the slice 
$\{A\}\times\{\Phi\}\times\GGG_G$. Then we prove that 
the functional being proper along $\{A\}\times\{\Phi\}\times\GGG_G$ is 
equivalent to its having a critical point in 
$\{A\}\times\{\Phi\}\times\GGG_G$, thus achieving the proof of
Theorem \ref{main}

\section{The integral of the moment map}
\label{integral}

In this section we consider the following general situation.
Let $H$ be a Lie group which acts on a Kaehler manifold
$M$ respecting the Kaehler structure, and assume that
there exists a moment map $$\mu:M\to\hlie^*,$$ where
$\hlie=\Lie(H)$. Suppose that there exists the complexification
$L=H^{\CC}$ of $H$, and that the inclusion $\iota:H\to L$
induces a surjection $\iota_*:\pi_1(H)\exh \pi_1(L)$. 
Under this assumptions, we construct a functional
$$\Psi:M\times L\to\RR$$
which we call the integral of the moment map $\mu$, and 
which satisfies these two properties:
\begin{itemize}
\item for any $x\in M$, the critical points of the restriction
$\Psi_x$ of $\Psi$ to $\{x\}\times L$ coincide with the points
of the orbit $Lx$ on which the moment map vanishes and
\item the restriction of $\Psi_x$ to lines of the form $\{e^{ts}|t\in\RR\}$, 
where $s\in\llie=\Lie(L)$, is convex.
\end{itemize}
If $H$ is compact then $L=H^{\CC}$ always exists and $\pi_1(H)\exh\pi_1(L)$
is always satisfied. But note that we do not need our manifold $M$ or our 
groups $H,\ L$ to be finite dimensional. In fact, we will use this 
construction mainly in the infinite dimensional case
$(M;H,L)=(\AAA^{1,1}\times\SSS;\GGG_K,\GGG_G)$ (in section 
\ref{esKaehler} we will prove that $\AAA^{1,1}\times\SSS$
is a Kaehler manifold, that the action of $\GGG_K$ respects the
Kaehler structure, and we will identify a moment map for this
action). The resulting integral of the moment map will be  
a certain modification of Donaldson functional.

\subsection{Definition of $\Psi$}

Let us fix a point $x\in M$, and let $\phi:L\to M$ be the map which sends
$h\in L$ to $hx\in M$. We define a 1-form on $L$,
$\sigma=\sigma^x\in\Omega^1(L)$, as follows: given
$h\in L$ and $v\in T_hL$, 
$$\sigma_h(v)=\la\mu(hx),-\imag \pi(v)\ra_{\klie},$$
where $\pi:T_hL=\hlie \oplus\imag \hlie \to\imag \hlie $
is the projection to the second summand.

We will use the following formula, which holds for any two 
vector fields $X,Y$ and any 2-form $\omega$ on $M$ 
\begin{equation}
d\omega(X,Y)=L_X(\omega(Y))-L_Y(\omega(X))-\omega([X,Y]).
\label{der}
\end{equation}
Equality (\ref{der}) is a particular case of a formula which
describes the exterior derivative of forms of arbitrary degree
in terms of Lie derivatives (see \cite{BeGeV} p. 18).

\begin{lemma}
The 1-form $\sigma$ is exact.
\end{lemma}
\begin{pf}
Let us first of all prove that $d\sigma=0$. Given $g\in\llie$, 
let $\fX^L_g\in\Gamma(TL)$ be the field generated
by $g$ acting on the left on $L$ (on the other hand, $\fX_g$
will denote the vector field generated by $g$ on $M$). 
We will prove that for any pair
$g,g'\in\hlie\cup\imag \hlie$, $d\sigma(g,g')=0$. This implies
by linearity that $d\sigma=0$. We will treat separately three cases,
and will make use of formula (\ref{der}), which in our case reads
$$d\sigma(\fX^L_g,\fX^L_{g'})=\la d(\sigma(\fX^L_{g'})),\fX^L_g\ra_{TL}
-\la d(\sigma(\fX^L_g)),\fX^L_{g'}\ra_{TL}-\sigma([\fX^L_g,\fX^L_{g'}]).$$

Suppose first that $g,g'\in\hlie$. In this case,
$\pi(\fX^L_g)=\pi(\fX^L_{g'})=\pi([\fX^L_g,\fX^L_{g'}])=0$, hence 
by the formula it is clear that $d\sigma(\fX^L_g,\fX^L_{g'})=0$.

Now suppose that $g\in\hlie$ and $g'\in\imag \hlie$. Observe that
$\sigma(\fX^L_g)=0$, so we have to prove that 
$\la d(\sigma(\fX^L_{g'})),\fX^L_g\ra_{TL}-\sigma([\fX^L_g,\fX^L_{g'}])=0$.
Differentiating property (C2) of the moment map (see section \ref{intro})
we have $\la d\la\mu,v\ra_{\hlie},\fX_g\ra_{TM}+\la\mu,[g,v]\ra_{\hlie}=0.$
The functoriality of the differentiation $d$ implies that 
$\la d(\sigma(\fX^L_{g'})),\fX^L_g\ra_{TL}+\sigma(\fX^L_{[g,g']})=0.$
On the other hand, since the action of $L$ on $M$ is on the left,
$[\fX^L_g,\fX^L_{g'}]=-\fX^L_{[g,g']}$
(see for example \cite{BeGeV} p. 208), so we 
obtain $$\la d(\sigma(\fX^L_{g'})),\fX^L_g\ra_{TL}-
\sigma([\fX^L_g,\fX^L_{g'}])=0,$$
which is what we wanted to prove.
The case $g\in\imag \hlie$ and $g'\in\hlie$ is dealt with in
a very similar way.

Finally, there remains the case $g,g'\in\imag \hlie$.
In this situation $[g,g']\in\hlie$, and so $\sigma([\fX^L_g,\fX^L_{g'}])=0$.
In view of this we have to prove 
$$\la d(\sigma(\fX^L_{g'})),\fX^L_g\ra_{TL}=
\la d(\sigma(\fX^L_g)),\fX^L_{g'}\ra_{TL}.$$
The left hand side is equal to 
$\phi^*(\la d\la\mu,\imag g\ra_{\hlie},\fX_{g'}
\ra_{TM})$ and this, by property (C1) of the moment map, is equal
to $\phi^*(\omega_M(I\fX_g,\fX_{g'}))=\phi^*(-\la \fX_g,\fX_{g'}\ra),$
where $\omega_M$ denotes the symplectic form on $M$.
The right hand side is equal to 
$\phi^*(\omega_M(I\fX_{g'},\fX_g))=\phi^*(-\la \fX_{g'},\fX_g\ra).$
Both functions are the same by the symmetry of $\la,\ra$.

Once we know that $d\sigma=0$, let us prove that $\sigma$ is exact.
Let $\iota:H\to L$ denote the inclusion. It is clear
that $\iota^*\sigma=0$. On the other hand, by our hypothesis
$\iota_*:\pi_1(H)\to\pi_1(L)$
is exhaustive. 
These two facts imply that $\sigma$ is exact.
Indeed, if it were not exact then we could find a path 
$\gamma:[0,1]\to L$, $\gamma(0)=\gamma(1)=1\in L$ such that 
$$\int_{\gamma}\sigma\neq 0.$$
But then we could deform $\gamma$ to a path $\gamma'\subset H$,
and, since $d\sigma=0$, the value of the integral would not change
and in particular would be nonzero. This is in contradiction with
the fact that $\iota^*\sigma=0$. So $\sigma$ is exact.
\end{pf}

Let $\Psi_x:L\to\RR$ be the unique function such that
$\Psi_x(1)=0$ and such that $d\Psi_x=\sigma^x$.
Define also $\Psi: M\times L \ni (x,g) \mapsto \Psi_x(g)$.
We will call the function $\Psi$ the {\rm integral of the moment map}.

\begin{remark}
If the symplectic form in $M$ is the curvature of a line bundle $L\to M$
and there is a lift of the action of $G$ to $L$, then the integral $\Psi$
of the moment map coincides with the functional defined in section
6.5.2 of \cite{DoKr}.
\end{remark}

\subsection{Properties of $\Psi$}
\label{propietats}
In this subsection we give the properties of the integral of the
moment map which will be used below.

\begin{prop} Let $x\in M$ be any point, and let $s\in\hlie$. 
\begin{enumerate}
\item $\Psi(x,e^{\imag s})=
\int_0^1 \la\mu(e^{\imag ts}x),s\ra_{\hlie} dt=
\int_0^1 \lambda_t(x;s) dt,$
\item $\frac{\partial\Psi}{\partial t}(x,e^{\imag ts})|_{t=0}
=\la\mu(x),s\ra_{\hlie}=\lambda_0(x;s),$
\item $\forall t_0\in\RR$,
$\frac{\partial^2\Psi}{\partial t^2}(x,e^{\imag ts})|_{t=t_0}\geq 0,$
with equality if and only if $\fX_s(e^{\imag t_0s}x)=0$,
\item $\forall t_0>0$, 
$\Psi(x,e^{\imag ls}x)\geq (l-t_0)\lambda_{t}(x;s)+C_s(x;t_0)$,
where $C_s(x;t_0)$ is a continuous function on $x\in M$,
$s\in\hlie$ and $t_0\in\RR$,
\end{enumerate}
\label{propietats1}
\end{prop}
\begin{pf}
By definition, $\Psi(x,e^{\imag s})=\int_{\gamma}\sigma^x$, where
$\gamma$ is any path in $L$ joining $1\in L$ to $e^{\imag s}$.
If we take $\gamma:[0,1]\ni t\mapsto e^{\imag ts}$, then
the integral reduces to $\int_0^1 \la\mu(e^{\imag ts}x),s\ra_{\hlie} dt$.
This proves (1). 
Property (2) is deduced from (1) differentiating. 
(3) is a consequence of (1) and the fact that $\lambda_t(x;s)$
increases with $t$.
To prove (4), let $C_s(x;t_0)=\int_0^{t_0}\lambda_t(x;s)dt$. Then:
$$\int_0^1 \lambda_t(x;ls)dt=\int_0^l \lambda_t(x;s)dt
\geq (l-t_0)\lambda_{t}(x;s)+C_s(x;t_0);$$
the first equality is obtained making a change of variable and 
using (2) in Proposition \ref{proppes}, and the 
inequality comes from the fact that $\lambda_t(x;s)$ increases as 
a function of $t$. 
\end{pf}

\begin{prop} Let $x\in M$ be any point, and let $s\in\hlie$. 
\begin{enumerate}
\item If $g,h\in L$, then
$\Psi(x,g)+\Psi(gx,h)=\Psi(x,hg)$,
\item for any $k\in H$ and $g\in L$, $\Psi(x,kg)=\Psi(x,g)$,
and $\Psi(x,1)=0$,
\item for any $k\in H$ and $g\in L$, $\Psi(kx,h)=\Psi(x,k^{-1}gk)$.
\end{enumerate}
\label{propietats2}
\end{prop}
\begin{pf}
To prove (1), observe that for any $g\in L$, $\sigma^{gx}=R_g^*\sigma^x$,
where $R_g$ denotes right multiplication in $L$ (indeed, 
for any $g'\in L$ one has $\sigma^{gx}(g')=\sigma^{x}(g'g)$
-- as usual, we identify the tangent spaces $T_{g'}(L)$
and $T_{g'g}(L)$ making $L$ act on the right).
This equivalence, together with the requierement that 
$\Psi_{gx}(1)=0$ implies that, for any $h\in L$, 
$\Psi_{gx}(h)=\Psi_x(hg)-\Psi_x(g)$. 
Property (2) is a consequence of (1) together with the fact that, for any 
$x\in M$, $\Psi_x|_H=0$. Finally, to prove (3) we use points (1) and (2):
$\Psi(x;k^{-1}gk)=\Psi(x,gk)+\Psi(gkx,k^{-1})
=\Psi(x,k)+\Psi(kx,g)=\Psi(kx,g).$
\end{pf}

\begin{prop} 
An element $g\in L$ is a critical point of $\Psi_x$ if and only 
if $\mu(gx)=0$.
\label{propietats3}
\end{prop}
\begin{pf} This is a consequence of (2) in \ref{propietats1}
and (1) in \ref{propietats2}.
\end{pf}

Just like maximal weights, the function $\Psi$ depends on the
moment map, which is not unique. When it is not clear from
the context which moment map we consider, we will write
$\Psi^{\mu}$ to mean the integral of the moment map $\mu$.

\subsection{Linear properness}
In this section we restrict to the case
$(M;H,L)=(F;K,G)$. 
Let $\rho_a:\glie\to\End(W_a)$ be the complexification of the (differential
of the) representation $\rho_a:K\to U(W_a)$ chosen in subsection 
\ref{isotriat}. We define a norm on $\glie$ as follows: for any $s\in\glie$,
$$|s|=\la s,s\ra^{1/2}=\Tr(\rho_a(s)\rho_a(s)^*)^{1/2}.$$
Let $\log_G:G\simeq K\times\exp(\imag\klie)\to\imag\klie$
denote the projection to the second factor of the Cartan decomposition
composed with the logarithm. For any $g\in G$ we will call 
$|g|_{\log}:=|\log_G g|$ the {\rm length} of $g$.

\begin{definition}
We will say that $\Psi_x$ is {\rm linearly proper} if there exist
positive constants $C_1$ and $C_2$ such that for any $g\in G$    
$$|g|_{\log}\leq C_1\Psi_x(g)+C_2.$$
\end{definition}

\begin{prop}
Let $h\in G$ and $x\in F$. If $\Psi_x$ is linearly proper 
then $\Psi_{hx}$ is also linearly proper.
\label{tambe}
\end{prop}

Before giving the proof of this proposition we prove the following
technical result.
\begin{lemma} Let 
$N=\dim W_a$ and $h\in G$. There exists $C\geq 1$ such that
for any $g\in G$
$$N^{-1/2}|gh|_{\log}-\log C
\leq |g|_{\log}\leq N^{1/2}(|gh|_{\log}+\log C).$$
Furthermore, $C$ depends continuously on $h\in G$.
\label{comparison}
\end{lemma}
\begin{pf}
Since the Cartan decomposition commutes with unitary representations,
we may describe the length function as follows. Let $x\in G$ be any element
and write $\rho_a(x)=RS$, where $R\in U(W_a)$ and $S=\exp(u)$,
where $u=u^*$. The matrix $u$ diagonalises and has real
eigenvalues $\lambda_1,\dots,\lambda_N$. So
$|x|_{\log}^2=\sum_{j=1}^N \lambda_j^2.$
Define $\max (x)=\max_{\|v\|=1}|\log \|\rho_a(x)v\||$. Then
we have $\max |\lambda_j|=\max (x)$ and consequently
\begin{equation}
\max(x)\leq |x|_{\log}\leq N^{1/2}\max(x).
\label{ee1}
\end{equation}
Let now $h\in G$. Then there exists $C\geq 1$, depending continuously
on $h$, such that for any $g\in G$ and any $v\in V$,
$C^{-1}\|\rho_a(gh)v\|\leq \|\rho_a(g)v\|\leq 
C\|\rho_a(gh)v\|$, which implies
\begin{equation}
|\max(gh)-\max(g)|\leq \log C.
\label{ee2}
\end{equation}
Putting $x=gh$ in (\ref{ee1}) we obtain
\begin{equation}
N^{-1/2}|gh|_{\log}\leq\max(gh)\leq|gh|_{\log},
\label{ee3}
\end{equation}
and combining (\ref{ee1}) with $x=g$ and (\ref{ee2}) we get
$$\max(gh)-\log C\leq |g|_{\log}\leq N^{1/2}(\max(gh)+\log C).$$
Finally, using (\ref{ee3}) we get
$N^{-1/2}|gh|_{\log}-\log C
\leq |g|_{\log}\leq N^{1/2}(|gh|_{\log}+\log C).$
\end{pf}

\begin{pf} (Proposition \ref{tambe})
Suppose that $\Psi_x$ is linearly proper, that is, for any $g\in G$
$$|g|_{\log}\leq C_1\Psi_x(g)+C_2,$$ where $C_1,C_2$ are positive.
Fix $h\in G$. Let $C\geq 1$ be the constant in Lemma \ref{comparison}.
(1) in \ref{propietats2} tells us that
$\Psi_{hx}(g)=\Psi_x(gh)-\Psi_x(h)$, so we get for any $g\in G$
\begin{align}
|g|_{\log} &\leq N^{1/2}(|gh|_{\log}+\log C)
\leq N^{1/2}(C_1\Psi_x(gh)+C_2+\log C)\notag\\
&=N^{1/2}(C_1(\Psi_x(gh)-\Psi_x(h))+C_1\Psi_x(h)+C_2+\log C)\notag\\
&=N^{1/2}(C_1\Psi_{hx}(g)+C_1\Psi_x(h)+C_2+\log C),\notag
\end{align}
so setting 
$C_1'=N^{1/2}C_1$ and $C_2'=\max\{0,N^{1/2}(C_1\Psi_x(h)+C_2+\log C)\}$ 
then $C_1',C_2'$ are positive and $|g|_{\log}\leq C_1'\Psi_{hx}(g)+C_2'.$
This proves that $\Psi_{hx}$ is linearly proper.
\end{pf}

\section{A Kaehler structure on $\AAA^{1,1}\times\SSS$}
\label{esKaehler}
In this section we will give, following the classical idea of Atiyah and  
Bott \cite{AB}, a $\GGG_K$-invariant Kaehler structure on 
the manifold $\AAA\times\SSS$. 
This structure will depend on our choice of a biinvariant metric on $\klie^*$,
and consequently on the representation $\rho_a$ used to define it.
We will identify for this structure a moment map of the action of $\GGG_K$, 
the maximal weights and the integral of the moment map. 

\subsection{Unitary connections}
\label{exemp1}

\subsubsection{$\AAA$ is a Kaehler manifold}
Let $\AAA$ be the space of 
$K$-connections on $E$. It is an affine space modelled on 
$\Omega^1(E\times_{\Ad}\klie)$.
We define a complex structure $I_{\AAA}$ on $\AAA$ as follows.
Given any $A\in\AAA$, the tangent space $T_A\AAA$ can be 
canonically identified with 
$\Omega^1(E\times_{\Ad}\klie)=\Omega^0(T^*X\otimes E\times_{\Ad}\klie)$.
Then we set $I_{\AAA}=-I^*\otimes 1$. The complex structure
$I_{\AAA}$ is integrable.
We also define on $\AAA$ a symplectic form $\omega_{\AAA}$.
Let $\Lambda:\Omega^{p,q}(X)\to\Omega^{p-1,q-1}(X)$ be the adjoint of
the map given by wedging with $\omega$. Then, if $A\in\AAA$ and 
$\alpha,\beta\in T_A\AAA\simeq\Omega^1(E\times_{\Ad}\klie)$, we set
$$\omega_{\AAA}(\alpha,\beta)=\int_X \Lambda(B_1(\alpha,\beta)).$$

Here $B_1:\Omega^1(E\times_{\Ad}\klie)\otimes\Omega^1(E\times_{\Ad}\klie)
\to \Omega^2$ is the combination of the usual wedge product with
the biinvariant nondegenerate pairing $\la,\ra$ on $\klie$ obtained
from the representation $\rho_a$.
It turns out that $\omega_{\AAA}$ is a symplectic form on $\AAA$, and it
is compatible with the complex structure $I_{\AAA}$. Hence $\AAA$ is
a Kaehler manifold. Furthermore, the action of $\GGG_G$ on 
$\AAA$ defined in subsection \ref{setting} is holomorphic and is the 
complexification of the action of $\GGG_K$.

\subsubsection{The moment map}
Recall that the Lie algebra of $\GGG_K$ is
$\Lie(\GGG_K)=\Omega^0(E\times_{\Ad}\klie)$.
There exists a moment map for the action of $\GGG_K$ on $\AAA$, and
it takes the following form (see for example \cite{DoKr, Ko}):
$$\begin{array}{rcl}
\mu:\AAA & \longrightarrow & \Lie(\GGG_K)^* \\
A & \mapsto & \Lambda F_A.
\end{array}$$
The curvature $F_A$ of $A$ lies in $\Omega^2(E\times_{\Ad}\klie)$, so 
$\Lambda F_A\in\Omega^0(E\times_{\Ad}\klie)\subset
\Omega^0(E\times_{\Ad}\klie)^*$, the last inclusion being
given by the integral of the pairing $\la,\ra$.

The proof of the next lemma is an easy exercise.
\begin{lemma}
Let $A\in\AAA$ be a connection, and take 
$s\in\Lie(\GGG_K)=\Omega^0(E\times_{\Ad}\klie)$. 
Then 
\begin{equation}
\lambda_t(A;s)=
\int_X\la\Lambda F_A,s\ra + 
\int_0^t \|e^{\imag ls}\overline{\partial}_A(s)e^{-\imag ls}\|^2 dl.
\end{equation}
\label{tpesmaxconn}
\end{lemma}
When $s\in L^2_1(E\times_{\Ad}\klie)$ the maximal weight is given
by exactly the same formula. To prove it
one needs to use a technical theorem of Uhlenbeck and Yau \cite{UY}. 
This result allows to regard $s$ as a genuine smooth section 
of $E\times_{\Ad}\klie$ at the complementary of a complex codimension
two subvariety of $X$, and to check that the integrals appearing
in Lemma \ref{tpesmaxconn} converge.

\subsubsection{The integral of the moment map}

The $K$-equivariance of the Cartan decomposition
implies that $\GGG_G\simeq \GGG_K\times \imag\Lie(\GGG_K),$
and from this fact, using that $\pi_1(K)\to\pi_1(G)$ is surjective,  
we see that $\pi_1(\GGG_K)\to\pi_1(\GGG_G)$ is a surjection (both
maps are the ones induced by the inclusions).
As a consequence, the results of section \ref{integral} apply to 
actions of $\GGG_K$ on Kaehler manifolds.
So there is an integral of the moment map 
$\Psi^{\AAA}$ which satisfies all the properties given in section
\ref{propietats}. Fix now a connection $A\in\AAA$. By (\ref{tpesmaxconn}) 
and using (1) in Proposition \ref{propietats1} 
we see that
\begin{align}
\Psi^{\AAA}_A(e^{\imag s})&=\int_0^1\lambda_t(A,s)
=\int_X\la\Lambda F_A,s\ra+                             
\int_0^1\left(
\int_0^t \|e^{\imag ls}\overline{\partial}_A(s)e^{-\imag ls}\|^2 dl 
\right)dt\notag\\
&=\int_X\la\Lambda F_A,s\ra+\int_0^1
(1-l)\|e^{\imag ls}\overline{\partial}_A(s)e^{-\imag ls}\|^2 dl.
\end{align}

Then, by (2) in \ref{propietats2}, the function $\Psi^{\AAA}_A$ 
factors through $$\Psi^{\AAA}_A:\GGG_G/\GGG_K\to\RR.$$
The resulting functional may be seen as a {\it modified
Donaldson functional}. In fact, when $F=\{\pt\}$,
it coincides (up to a multiplicative constant) with the
Donaldson functional. To see this, one only has to check 
that the Donaldson functional satisfies property (2) in 
\ref{propietats1} (see \cite[Lemma 3.3.2]{Br2} for the case
$F=\CC^n$).

We will use the restriction of the integral of the moment map to
$\AAA^{1,1}\times\SSS\times\GGG_G$ ($\AAA^{1,1}\subset\AAA$ is a 
conplex subvariety, but in general it is not smooth).
This functional will be the main tool in proving Theorem \ref{main}.

\subsubsection{Maximal weights for $A\in\AAA^{1,1}$}
\label{pesmaxconn2}
Note that since $\AAA^{1,1}\subset\AAA$ is a $\GGG_G$ invariant
subvariety, the moment map, the maximal weights and the integral
of the moment map of the action of $\GGG_K$ on $\AAA^{1,1}$
are the restrictions of their counterparts in $\AAA$.

Recall that $V=E\times_{\rho_a}W_a\to X$ is the vector bundle associated to
the representation $\rho_a$. For any $s\in\Lie(\GGG_K)$ we can
view $\rho_a(s)$ as a section of $E\times_{\Ad(\rho_a)}\End(W_a)$.
Take a connection $A\in\AAA^{1,1}$, and consider on $V$ the holomorphic
structure induced by $\overline{\partial}_A$. 
Using Lemma \ref{tpesmaxconn} one can prove the following (see \cite{Mu}).
\begin{lemma}
Let $s\in\Lie(\GGG_K)$. If $\lambda(A;s)<\infty$, then 
the eigenvalues of $\rho_a(s)$ are constant.
Let $\lambda_1<\dots<\lambda_r$ be the different eigenvalues
of $\imag\rho_a(s)$,
and let $V(\lambda_j)\subset V$ be the eigenbundle of eigenvalue
$\lambda_j$. Put 
$V^{\lambda_k}=\bigoplus_{j\leq k} V(\lambda_j)$.
Then, for any $k$, $V^{\lambda_k}$ is a holomorphic subbundle
of $V$. Furthermore
$$\lambda(A;s)=
\lambda_r\deg(V)+\sum_{k=1}^{r-1}(\lambda_k-\lambda_{k+1})
\deg(V^{\lambda_k}).$$
\label{redhol}
\end{lemma}

If we consider more generally $s\in L^2_1(E\times_{\Ad}\klie)$,
then $\lambda(A;s)<\infty$ leads to a filtration of the locally
free sheaf associated to $V$ by reflexive (coherent) subsheaves,
and not only holomorphic subbundles of $V$ as in the smooth
case. To prove this one uses a theorem of Uhlenbeck and Yau (see \cite{UY}
and \cite[\S 3.11]{Br2}).

\subsection{Sections of the associated bundle}
\label{exemp2}

\subsubsection{$\SSS$ is a Kaehler manifold}
Let us define a complex structure $I_{\SSS}$ and a symplectic
form $\omega_{\SSS}$ on $\SSS=\Gamma(\FFF)$.
Consider a section $\sigma\in\SSS$. The tangent space 
$T_{\sigma}\SSS=\Gamma(\sigma^* T\FFF_v)$,
where $T\FFF_v\subset T\FFF$ is the subbundle of vertical tangent vectors
of $\FFF$, that is, $T\FFF_v=\Ker(d\pi_F)$. 
Let $\alpha\in\Gamma(\sigma^* T\FFF_v)$. We set by definition
$I_{\SSS}(\alpha)=I_F\alpha$. This makes sense, since the $K$ invariance
of $I_F$ implies that $T\FFF_v$ inherits the complex structure of $F$.
Now let $\alpha,\beta\in\Gamma(\sigma^* T\FFF_v)$. We
define the symplectic form $\omega_{\SSS}$ on $\SSS$ as
$$\omega_{\SSS}(\alpha,\beta)=\int_X \omega_F(\alpha,\beta).$$
The 2-form $\omega_{\SSS}$ is nondegenerate (this 
is a consequence of the nondegeneracy of $\omega_F$) and 
$\omega_{\SSS}$ and $I_{\SSS}$ are compatible, that is,
$\la\alpha,\beta\ra=\omega_{\SSS}(\alpha,I_{\SSS}\beta)$
is a Riemannian pairing. The two structures are integrable, and
so $\SSS$ is a Kaehler manifold.

\subsubsection{The actions of $\GGG_K$ and $\GGG_G$ and the moment map} 
Both groups $\GGG_K$ and $\GGG_G$ act on the space
of sections $\SSS=\Gamma(\FFF)$, and the action of
$\GGG_G$ is the complexification of the action of $\GGG_K$.
On the other hand, $\GGG_K$ acts by isometries and respecting the 
symplectic form, and there exists a moment map $\mu_{\SSS}$, which 
is equal fibrewise to $\mu$ (the moment map of the action of $K$ on $F$). 
As such, it is a section of $\Omega^0(E\times_{\Ad}\klie)^*$.

\subsubsection{Maximal weights}
\label{pesmaxsec}
The maximal weight of $s\in\Lie(\GGG_K)=\Omega^0(E\times_{\Ad}\klie)$
acting on a section $\Phi\in\SSS$ is given by the integral of
the maximal weight in each fibre:
$$\int_{x\in X}\lambda(\Phi(x);s(x)).$$
This makes sense due to the $K$ equivariance of $\lambda$. See
(1) in Lemma \ref{proppes}.

\subsubsection{The integral of the moment map} The results
in section \ref{integral} imply that there exists an integral
$\Psi^{\SSS}$ of the moment map of the action of $\GGG_K$ on $\SSS$. If 
$\Psi:F\times G\to\RR$ is the integral of the moment map of the
action of $K$ on $F$, then, for any section $\sigma\in\SSS$ and
gauge transformation $g\in\GGG_G$
$$\Psi^{\SSS}(\sigma,g)=\int_{x\in X} \Psi(\sigma(x),g(x)).$$
This makes sense due to the $K$-equivariance of $\Psi$:
see (3) in \ref{propietats2}. 

\subsection{Symplectic point of view}
\label{disgressio}

We saw that both $\AAA^{1,1}$ and 
and $\SSS$ are Kaehler manifolds, with symplectic forms $\omega_{\AAA}$
and $\omega_{\SSS}$ and with actions of $\GGG_K$ extending
to actions of the complexification $\GGG_G$. Hence
$\AAA^{1,1}\times\SSS$ is also a Kaehler manifold, with symplectic
form $\omega_{\AAA}+\omega_{\SSS}$ (we omit the pullbacks). The moment map
$\mu_{\AAA\times\SSS}$
of the action of $\GGG_K$ on $\AAA\times\SSS$ will simply be the moment
map of the action on $\AAA$ plus that of the action on $\SSS$.
That is, $$\mu_{\AAA\times\SSS}(A,\Phi)=\Lambda F_A+\mu(\Phi).$$
So equation (\ref{equacio}) can be written as 
$\mu_{\AAA\times\SSS}=c$, where $c$ denotes the central element
in $(\Lie(\GGG_K))^*=\Omega^0(E\times_{\Ad}\klie)^*$ which is
fibrewise equal to a central element $c\in\klie^*$. 
Furthermore, we have the following result.
\begin{lemma} $T^c_{\Phi}(\sigma,\chi)=
\lambda^{\Lambda F_A+\mu(\Phi)-c}((A,\Phi);-\imag g_{\sigma,\chi}).$
\label{equivalencia}
\end{lemma}
\begin{pf}
Combine subsections \ref{pesmaxconn2} and \ref{pesmaxsec}.
\end{pf}

As a final comment, note that so far we have defined the gauge 
group as the space of smooth sections of a certain bundle. Eventually, 
it will be necessary to take a metric on $\GGG_K$ (and $\GGG_G$) and 
complete both spaces with respect to the metric, to assure the convergence
of certain sequences. We will use Sobolev $L^p_2$ and $L^2_1$ norms.

\section{Analytic stability and vanishing of the moment map in finite 
dimension}
\label{correspondencia}
We will now pause to prove Theorem \ref{main} in the case
$X=\{\pt\}$, which is much easier than the general case and is 
interesting {\it per se}.
The results in this section (at least for the case in which 
$F$ is projective) have been known for many years: see \cite{KeNe, Ki}.
That they are related with Hitchin--Kobayashi correspondence
was also known since the first cases of the correspondence were
studied. Our intention here is to make more concrete this relation and 
to stress on the similarities between the
{\it finite dimensional situation} $X=\{\pt\}$ and the general one considered
in Theorem \ref{main} (which corresponds to the situation in which
$F=\AAA^{1,1}\times\SSS$ with the actions of $\GGG_K$ and $\GGG_G$). 
For example, the different versions of
Donaldson functional used in the literature are in fact particular
instances of a construction which works for a wide class of Kaehler
actions of Lie groups on Kaehler manifolds (namely, what we have
called the integral of the moment map). Moreover, the $c$-stability
condition is also a particular case of a general notion of stability
for group actions on Kaehler manifolds (the so-called analytic
stability). And the very correspondence coincides almost word by word
with Theorem \ref{corr} given in this section. The proof which we
give here works only for Kaehler actions of compact groups,
and so it can not be used in the general situation (in which the group
is $\GGG_K$). Nevertheless, the scheme of the proof will be the same
in the general situation.

Let us write $\Psi:F\times G\to\RR$ for the integral of the moment
map $\mu:F\to\klie^*$.

\begin{definition}
Let $x\in F$. We will say that $x$ is {\rm analytically stable} 
if for any $s\in\klie$ the maximal weight of $s$ acting on $x$ 
is strictly positive: $$\lambda(x;s)>0.$$
\label{estabilitat}
\end{definition}

\begin{lemma} 
A point $x\in F$ is analytically stable if and only if 
$\Psi_x$ is linearly proper.
\label{fita}
\end{lemma}
\begin{pf}
Suppose first that $x$ is analytically stable.
We have to prove that there exists two positive constants
$C_1,C_2\in\RR$ such that,
for any $s\in\klie$, $|s|\leq C_1\Psi_x(e^{\imag s})+C_2$.
Assume that there are not such constants. Then, we can find sequences
$\{s_j\}\subset \klie$ and $\{C_j\}\subset\RR$ such that 
$|s_j|\to\infty$, $C_j\to\infty$ and, for any $j$, 
$|s_j|\geq C_j\Psi_x(e^{\imag s_j})$. Let $u_j=s_j/|s_j|$.
After passing to a subsequence, we can assume that 
$u_j\to s$. Take now any $t>0$. By our hypothesis, and making use
of (4) in Proposition \ref{propietats1}, 
$$\frac{1}{C_j}\geq\frac{\Psi_x(e^{\imag s_j})}{|s_j|}
\geq \frac{(|s_j|-t)}{|s_j|}
\lambda_{t}(x;u_j)+\frac{C_{u_j}(x;t)}{|s_j|}.$$
Now, making $j\to\infty$, we obtain $0\geq \lambda_{t}(x;s)$,
since, by the compactness of $B_{\klie}(1)=\{s\in\klie|\ |s|=1\}$,
$C_{u_j}(x;t)$ is uniformly bounded.
This is true for any $t>0$, so passing to the limit $t\to\infty$
we get $0\geq \lambda(x;s)$, which contradicts analytic stability.

Now suppose that there exists positive $C_1,C_2$ such that for any 
$s\in\klie$ 
\begin{equation}
|s|\leq C_1\Psi_x(e^{\imag s})+C_2.
\label{supo}
\end{equation}
We have to prove that $x$ is analytically stable.
So take $s\in\klie$ and assume that $\lambda(x;s)\leq 0$.
In this case, for any $t\geq 0$, 
$\Psi_x(e^{\imag ts})=\int_0^t\lambda_l(x;s)dl\leq 0$,
which, for $t$ big enough, contradicts (\ref{supo}).
This proves that $x$ is analytically stable.
\end{pf}

\begin{corollary}
Let $x\in F$. Then $x$ is analytically stable if and only if $hx$ is 
analytically stable for any $h\in G$.
\label{tambe2}
\end{corollary}
\begin{pf}
This is a consequence of the preceeding lemma together
with Lemma \ref{tambe}.
\end{pf}

\begin{theorem}
Let $x\in F$ be any point. There is at most one $K$ orbit inside
the orbit $Gx\subset F$ on which the moment map vanishes. 
Furthermore, $x$ is analytically stable if and only if:
(1) the stabiliser $G_x$ of $x$ in $G$ is finite and (2) there exists
a $K$ orbit inside $Gx$ on which the moment map vanishes.
\label{corr}
\end{theorem}
\begin{pf}
We first prove uniqueness.
Assume that there are two different $K$ orbits inside a $G$ orbit
on which the moment map vanishes, say $Kx$ and $Kgx$, where $g\in G$. 
By the polar decomposition we can assume that $g=e^{\imag s}$, 
where $s\in\klie$. Consider the function $\Psi_x:G\to\RR$. 
By Proposition \ref{propietats3}, since $\mu(x)=0$, 
both $1,g\in G$ are critical points of $\Psi_x$. Consider now
the path $\gamma(t)=e^{\imag ts}$ connecting $1$ and $g$.
(3) of the proposition tells us that the restriction $\psi$
of $\Psi_x$ to this path has second derivative $\geq 0$. 
Since $0$ and $1$ are critical points of $\psi$, 
the second derivative must vanish at any point between $0$ and $1$.
In particular, 
$\frac{\partial^2\Psi}{\partial t^2}(x,e^{\imag ts})|_{t=0}=0$;
but this implies (again, (3) of the proposition), that the vector field 
$\fX_s(x)=0$, which gives $\fX_{\imag s}(x)=I\fX_s(x)=0$. 
So $e^gx=e^{\imag s}x=x$, and the two orbits $Kgx$ and $Kx$ coincide.

Suppose now that the point $x$ is analytically stable. 
Let us see that there is a $K$ orbit inside $Gx$ on which $\mu$ vanishes. 
By Lemma \ref{fita}, the function $\Psi_x$ is linearly proper.
Using (2) in \ref{propietats2}, we
conclude that there must exist a critical point in the $G$
orbit of $x$. Indeed, if $\{s_j\}\subset\klie$ are such that
$e^{\imag s_j}$ is a minimising sequence for $\Psi_x$, then by
the preceeding lemma the set $\{s_j\}$ is bounded; so it
has a subsequence converging to a certain $s\in\klie$, and
$e^{\imag s}$ is a minimum of $\Psi_x$ (of course, here we
use that $\klie$ has finite dimension). At this point (even 
more, at the $K$ orbit through this point) the moment map must vanish.
Let now $y=e^{\imag s}x$. By Lemma \ref{tambe2} $y$ is analitically stable.
If the stabiliser $K_y$ of $y$ in $K$ were not finite, 
then, since $K$ is compact, its closure would be a Lie 
subgroup of $K$ of dimension greater than zero. In particular, there 
would exist an $s\in\klie$ such that $\fX_s(y)=0$. But then $e^{ts}y=y$ 
for any $t$, so that the gradient flow $\phi^t_s$ leaves $y$ fixed.
This means that $\lambda(y;s)=-\lambda(y;-s)$, so that either $\lambda(y;s)$
or $\lambda(y;-s)$ (or both) is $\leq 0$. This contradicts 
analytic stability. So $K_y$ is finite. 

Finally, since $\mu(y)$ is invariant under the coadjoint action of 
$K$ in $\klie^*$, it turns out that $G_y$ is the complexification of $K_y$. 
Let us see why (we copy the proof of \cite[Proposition 1.6]{Sj}).
One inclusion is easy: $G_y$ contains the complexification of $K_y$.
For the other inclusion,
let $ge^{\imag s}$ be an arbitrary element of $G_x$, where
$g\in K$ and $s\in\klie$. We want to show that $g\in K_x$ and
$s\in\klie_x$ (where $\klie_x$ is the infinitesimal stabliser
of $x$). Using the fact that $\mu$ is $K$-equivariant we have
$$\mu(e^{\imag s}x)=g^{-1}\mu(ge^{\imag s}x)=g^{-1}\mu(x)=\mu(x).$$
Now, Lemma \ref{gradient} implies that $s\in\klie_x$, from which we
deduce that $g\in K_x$. This finishes the proof.
So $G_y$ is finite and in consequence $G_x$ is also finite.

To prove the converse, let $x\in F$.
Assume that $G_x$ is finite and that there exists $g\in G$ such 
that $\mu(gx)=0$. Then $G_{gx}$ is also finite and consequently
so is $K_{gx}$. This implies that, for any $s\in\klie$,
$\fX_{\imag s}(gx)\neq 0$, so (Lemma \ref{gradient}),
$\lambda(gx;s)>\mu_s(gx)=0$. This means that $gx$ is analytically
stable, hence so is $x$.
\end{pf}

It is an exercise to verify that the property
on analitically stable points of $F$ of being simple 
(see subsection \ref{ssp}) is equivalent to that of having finite 
stabiliser in $G$.

Using the results in this section one can also study the equation
$\mu=c$, where $c\in\klie^*$ is any central element. 
Indeed, $\mu-c$ is a moment map, and so one only has to consider
the maximal weights $\lambda^{\mu-c}$ and the integral 
$\Psi^{\mu-c}$.

\begin{remark}
Suppose that $F\subset\PP^n$ is a projective manifold and 
that the Kaehler structure on $F$ is that induced by the Kaehler 
structure on $\PP^n$. Using the Hilbert-Mumford numerical criterion, 
one can easily prove that in this context the property of
being analytically stable and having finite stabiliser 
is the same as being stable in the sense of Mumford Geometric
Invariant Theory (see \cite{MFK} and lemma 8.8 and remark 8.9 in
\cite{Ki}).
\end{remark}

\section{Proof of the correspondence}
\label{demostracio}
\subsection{The length of elements of the gauge group}
There are several ways to extend the notion of length to elements of
the gauge group. We will use these two definitions: if
$g\in\GGG_G$, then $|g|_{{\log},C^0}=\| |g|_{\log}\|_{C^0}$ and
$|g|_{{\log},L^1}=\| |g|_{\log}\|_{L^1}$ (to give 
this a sense we use the $K$ invariance of the length function, which 
is a consequence of the fact that the Cartan decomposition 
$G\simeq K\times\exp(\imag\klie)$ is $K$-equivariant).
Define a norm $\|\cdot\|_{L^p}$ in 
$\Lie(\GGG_G)=\Omega^0(K\times_{\Ad}\glie)$ as 
the $L^p$ norm of $|\cdot|$: if 
$s\in\Omega^0(K\times_{\Ad}\glie)$ then
$$\|s\|_{L^p}=\left(\int_{x\in X} |s(x)|^p\right)^{1/p}.$$
We will usually write $\|\cdot\|$ instead of $\|\cdot\|_{L^2}$.

\subsection{Stability implies existence of solution}
Here we will follow the scheme in section \ref{correspondencia}.
Fix a pair $(A,\Phi)\in\AAA^{1,1}\times\SSS$.
We will make use of the integral of the moment map
$\mu^c(A,\Phi)=\Lambda F_A+\mu(\Phi)-c$,
$\Psi^c=(\Psi^{\AAA\times\SSS})^{\mu^c}_{(A,\Phi)}=
(\Psi^{\AAA})^{\mu^c}_A+(\Psi^{\SSS})^{\mu^c}_{\Phi}$, 
and will see that if the pair $(A,\Phi)$
is simple and $c$-stable, then there exists a $\GGG_K$ orbit inside
the $\GGG_G$ orbit of $(A,\Phi)$ on which $\Psi^c$ attains
its minimum. The main step will be to prove that if the
condition of $c$-stability is satisfied, then the map
$\Psi^c$ satisfies an inequality like that in Lemma \ref{fita}.
This method of proof is exactly the same that appears in
\cite{Si, Br2, BrGP1, DaUW} (and in many other places where
similar results are proved), though here we have tried to 
remark the similarities with the finite dimensional case,
so our notation changes a little bit. However, in some steps of the proof
we will only give a sketch, refering to \cite{Br2} for details. 

Recall that on $\glie$ we have a Hermitian
pairing $\la,\ra:\glie\times\glie\to\CC$ and a norm $|\cdot|$,
both obtained by means of the representation $\rho_a$.
We will use the following $L^p$ norm on $\Omega^0(E\times_{\Ad}\glie)$:
$$\|s\|_{L^p}=\left(\int_X |s(x)|^p\right)^{1/p},$$
and Sobolev norm
$$\|s\|_{L^p_2}=\|s\|_{L^p}+\|d_A s\|_{L^p}
+\|\nabla d_A s\|_{L^p},$$
where $\nabla:\Omega^0(T^*X\otimes E\times_{\Ad}\glie)\to 
\Omega^1(T^*X\otimes E\times_{\Ad}\glie)$ is
$\nabla_{LC}\otimes d_A$, $\nabla_{LC}$ being the Levi-Civita
connection. 
As usual, $L^p_2(E\times_{\Ad}\glie)$ will denote
the completion of $\Omega^0(E\times_{\Ad}\glie)$ with respect to
the norm $\|\cdot\|_{L^p_2}$.

\subsubsection{}
Suppose from now on that $(A,\Phi)$ is simple and $c$-stable. Our aim is to
minimise $\Psi^c$ in $\GGG_G/\GGG_K$. Through the exponential
map we can identify $\GGG_G/\GGG_K$ with 
$\Omega^0(E\times_{\Ad}\imag\klie)$.
Fix from now on $p>2n$ and define
$$\Met=L^p_2(E\times_{\Ad}\imag\klie).$$
The first thing to do is to restrict ourselves to the subset of 
$\Met$ defined as follows:
$$\MetB=\{s\in\Met\mid \|\mu^c(e^s(A,\Phi))\|^p_{L^p}\leq B\}.$$
Here $B$ is any strictly positive real constant. We prove that if a metric 
minimizes the functional in $\MetB$, then it also minimizes it in $\Met$.
For that it is enough to see that any minimum in $\MetB$
lies away from the boundary of $\MetB$; to verify this claim
one needs the hypothesis that the pair $(A,\Phi)$ is simple. Let us briefly 
explain how this goes (see also \cite{Br2}, Lemma 3.4.2).

Suppose that $s$ minimizes the functional inside $\MetB$. 
Let $B=e^s(A)$, $\Theta=e^s(\Phi)$. Define the differential operator
$L:L^p_2(E\times_{\Ad}\imag\klie)\to L^p(E\times_{\Ad}\imag\klie)$
as $$L(u)=\imag\left.\frac{\partial}{\partial t} 
\mu^c(e^{tu}(B,\Theta))\right|_{t=0}
=\imag\la d\mu^c,u\ra_{T(\AAA\times\SSS)}(B,\Theta).$$
Now, if we can see that there exists an $u$ such that
\begin{equation}
L(u)=-\imag\mu^c(B,\Theta), 
\label{lmh}
\end{equation}
then we can deduce that $\mu^c(B,\Theta)=0$ and, hence, that
$s$ minimizes the functional in the whole space of metrics
$\Met$ (see \cite[Lemma 3.4.2]{Br2} for a proof of this fact).
The operator $L$ is Fredholm and has index zero. Indeed,
modulo a compact operator it is 
$\imag\Lambda\overline{\partial}_{B}\partial_{B}$.
Using the Kaehler identities this is equal to
$\partial_{B}^*\partial_{B}$, which is an elliptic self adjoint operator.
This implies that if $\Ker(L)=0$ then $L$ is surjective 
and so, in particular, equation (\ref{lmh}) has a solution. 
Assume that $L(u)=0$, where $u\in\Met$. Then, by Lemma \ref{gradient}, 
\begin{align}
0&=\la -\imag L(u),-\imag u\ra=
\la \la d\mu^c,u\ra_{T(\AAA\times\SSS)},-\imag u\ra_{\Lie(\GGG_K)}
(B,\Theta)\notag\\
&=\|\fX^{\AAA^{1,1}\times\SSS}_{-\imag u}(B,\Theta)\|^2.
\end{align}
And this implies that $-\imag u$ leaves $(B,\Theta)$ invariant.
Hence if $u\neq 0$ then, since $u$ is semisimple, 
$(B,\Theta)$ is not simple, so neither is
$(A,\Phi)$; and this is a contradiction.

\subsubsection{}
The next step is to prove that the functional $\Psi^c$ is
linearly proper with respect to the $C^0$ norm in $\GGG_G$.

\begin{lemma}
There exist positive constants $C_1,C_2$ such that for any $s\in\MetB$ 
one has $\sup |s|\leq C_1\Psi^c(e^s)+C_2.$
\label{fita2}
\end{lemma}
\begin{remark}
It makes sense to speak about $\sup|s|$ because, since we took
$p>2n$, the Sobolev embedding theorem implies that 
$L^p_2\hookrightarrow C^0$ continuously (in fact, this is
a compact embedding).
\end{remark}
Just as in Lemma \ref{fita}, it is here that one uses the stability 
of the pair $(A,\Phi)$. First of all one sees that such a bound is 
equivalent to an $L^1$ bound: $\|s\|_{L^1}\leq C_1\Psi^c(e^s)+C_2$
(the constants in both inequalities need not be the same!). One uses that
pointwise
\begin{equation}
|s|\Delta|s|\leq \la\Lambda F_{e^s(A)} - \Lambda F_A,-\imag s\ra.
\label{equ1}
\end{equation}
This is proved in full detail in (\cite{Br2}, Prop. 3.7.1) for
$G=GL(n;\CC)$ and the metric induced by the fundamental representation. 
In our case, we use the representation $\rho_a$ to apply
this result to our $G$. 

\begin{lemma} For any point $x\in X$
\begin{equation}
0\leq \la\mu(e^s\Phi(x))-\mu(\Phi(x)),-\imag s(x)\ra_{\klie}.
\label{equ2}
\end{equation}
\end{lemma}
\begin{pf}
The gradient flow of $\mu_{-\imag s}$ is precisely $e^s$
(see Lemma \ref{gradient}).
\end{pf}

Summing the inequalities (\ref{equ1}) and (\ref{equ2}), 
using Cauchy-Schwartz, and dividing by $|s|$ we obtain the pointwise bound
$\Delta|s|\leq |\mu^c(e^s(A,\Phi))-\mu^c(A,\Phi)|$. 
Now, by of a result of Donaldson (see \cite{Br2}, Lemma 3.7.2), 
this bound allows to relate the $C^0$ and $L^1$ norms
of $s$ provided $s\in\MetB$. More precisely, we conclude that there
exists a constant $C_B$ such that for any $s\in\MetB$ one has
$\|s\|_{C^0}\leq C_B\|s\|_{L^1}$.

\subsubsection{}
In order to prove the existence of constants $C_1$ and $C_2$ such that
$\|s\|_{L^1}\leq C_1 \Psi^c(e^s)+C_2$, we suppose the contrary and try
to deduce that in this case the pair $(A,\Phi)$ cannot be $c$-stable.
If there exist not such constants, then we can find a sequence of
real numbers $C_j\to\infty$ and elements $s_j\in\MetB$ with
$\|s_j\|_{L^1}\to\infty$ such that 
$\|s_j\|_{L^1}\geq C_j \Psi^c(e^s)$
(see \cite{Br2}, Lemma 3.8.1). 
Set $l_j=\|s_j\|_{L^1}$, $u_j=l_j^{-1}s_j$ so that
$\|u_j\|_{L^1}=1$ and $\sup|u_j|\leq C$.

\begin{lemma}
After passing to a subsequence, there exists 
$u_{\infty}\in L^2_1(E\times_{\Ad}\imag\klie)$ such that
$u_j\to u_{\infty}$ weakly in $L^2_1(E\times_{\Ad}\imag\klie)$
and such that $\lambda((A,\Phi);-\imag u_{\infty})\leq 0.$
\label{contradic}
\end{lemma}
\begin{pf}
Just as in Lemma \ref{fita}, take $t>0$. Then (4) in proposition
\ref{propietats1} gives
\begin{align}
\frac{1}{C_j}&\geq\frac{\Psi^c(e^{s_j})}{\|s_j\|}\geq
\frac{l_j-t}{l_j}\lambda_t((A,\Phi);-\imag u_j)+\frac{1}{l_j}
\int_0^t\lambda_l((A,\Phi);-\imag u_j)dl\notag\\
&=\frac{l_j-t}{l_j}(\lambda_t(A;-\imag u_j)+\lambda_t(\Phi;-\imag u_j))
\notag\\
&+\frac{1}{l_j}\int_0^t(\lambda_l(A;-\imag u_j)+
\lambda_l(\Phi;-\imag u_j))dl.
\label{desclau}
\end{align}
Now, since $\|u_j\|_{C^0}\leq C_B$, 
and $X$ is compact, $\lambda_t(\Phi;-\imag u_j)$ and 
$\int_0^t\lambda_t(\Phi;-\imag u_j)dl$ are both bounded. Hence, there exists 
$C$ such that for any $j$ $$\frac{l_j-t}{l_j}\lambda_t(A;-\imag u_j)+
\frac{1}{l_j}\int_0^t\lambda_l(A;-\imag u_j)dl<C.$$
Using again the boundedness of $\|u_j\|_{C^0}$ and taking into
account Lemma \ref{tpesmaxconn} we obtain 
$$\|\overline{\partial}_A(u_j)\|_{L^2}<C_1.$$
Now, $\overline{u_j}=u_j$ (because the Cartan involution leaves
$\imag\klie$ fixed), and this implies that $\|u_j\|_{L^2_1}$ is also
bounded. So we can take a subsequence (which we again
call $\{u_j\}$) that converges weakly to $u_{\infty}\in L^2_1$. 
We can also assume that there exists the limit
$\lim_{i\to\infty} \lambda_t((A,\Phi);-\imag u_j)$.
On the other hand, since the embedding $L^2_1\hookrightarrow L^2$ 
is compact, we get strong convergence $u_j\to u_{\infty}$ in $L^2$.
$\|u_j\|_{L^1}=1$ and the uniform bound $\|u_j\|_{C^0}\leq C_B$
imply that $\|u_j\|_{L^2}>C_B^{-1}>0$, so $u_{\infty}\neq 0$.
To see that $\lambda_t((A,\Phi);-\imag u_{\infty})\leq
\lim_{i\to\infty} \lambda_t((A,\Phi);-\imag u_j)$ we observe that 
$$u_j\in L^2_{0,C_B}(E\times_{\Ad}\imag\klie)=
\{s\in L^2(E\times_{\Ad}\imag\klie)|\ |s(x)|\leq C_B\mbox{ a.e.}\}.$$
This implies that $u_{\infty}\in L^2_{0,C_B}(E\times_{\Ad}\imag\klie)$,
and this is enough to get the inequality (see \cite[Proposition 3.2.2]{Br2}). 
Finally, making $j\to\infty$ in formula (\ref{desclau}) we obtain
$$\lim_{i\to\infty} \lambda_t((A,\Phi);-\imag u_j)\leq 0,$$ so in particular
$\lambda_t((A,\Phi);-\imag u_{\infty})\leq 0$. Since this is true for any
$t>0$, we get $\lambda((A,\Phi);-\imag u_{\infty})\leq 0$.
\end{pf}

The next steps are rather standard. One can prove that
$\rho_a(u_{\infty})$ has almost everywhere constant eigenvalues and that it
defines a filtration of $V$ by holomorphic subbundles in the complement
of a complex codimension 2 subvariety of $X$. 
This follows exactly the same lines as \cite[\S\S 3.9, 3.10]{Br2},
the main technical point being the use of a theorem of Uhlenbeck and
Yau \cite{UY} on {\it weak subbundles} of vector bundles
(see \cite[\S 3.11]{Br2}).
The filtration of $V$ on $X_0$ and the gauge transformation 
$u_{\infty}$ lead to a reduction of the structure group
$\sigma\in\Gamma(X_0;E(G/P))$ defined on $X_0$ by \ref{filtsect}
which will be holomorphic thanks to the results in
subsection \ref{holomred}, and
an antidominant character $\chi$ of $P$. The degree of the pair 
$(\sigma,\chi)$ equals $\lambda((A,\Phi);-\imag u_{\infty})\leq 0$. 
And this contradicts the stability condition, thus finishing the proof
of Lemma \ref{fita2}.

\subsubsection{} With the inequality of Lemma \ref{fita2} in our hands,
we finish the proof of existence of solution to the equations
exactly as is done in \cite[\S 3.14]{Br2}. This consists
of two steps: the first one is to verify that there exists
an element $s\in\MetB$ minimising $\Psi^c$ and the second
one is to prove the smoothness of this solution $s$.

\subsection{Existence of solutions implies stability}

The method we will follow in this section will be exactly the
same as in the finite dimensional case in section \ref{correspondencia}.
Let us take a simple pair $(A,\Phi)\in\AAA^{1,1}\times\SSS$.
Suppose that there exists a gauge transformation $h\in\GGG_G$ 
such that $h(A,\Phi)$ satisfies equation (\ref{hk}).
We want to prove that $(A,\Phi)$ is analitically stable.

Take $X_0\subset X$ with complement of complex codimension 2, 
$P\subset G$ parabolic, $\chi$ an antidominant character of $P$ 
and fix a reduction $\sigma\in\Gamma(X_0;E(G/P))$. 
Thanks to \ref{sectfilt} we get a section 
$g_{\sigma,\chi}\in\Omega^0(X_0;E\times_{\Ad}\imag\klie)$,
and we have to check that $\lambda((A,\Phi);-\imag g_{\sigma,\chi})>0$.
In the following two lemmae it will be necessary to take into account
that $X_0$ has finite volume and that it has no nonconstant holomorphic 
functions (the last claim follows from Hartog theorem).

\begin{lemma}
For any semisimple $s\in L^p_2(X_0;E\times_{\Ad}\imag\klie)$ we have
$\lambda(h(A,\Phi);s)>0$.
\end{lemma}
\begin{pf}
Suppose the contrary: $\lambda(h(A,\Phi);s)\leq 0$.
Arguing as in \cite[\S 3.11]{Br2} (see also Lemma \ref{redhol})
we deduce that the eigenvalues of $s$ are constant.
Suppose that $s$ fixes $h(A,\Phi)$. Let $A'=h^*A$. Then
$d_{A'}s=0$, so $\overline{\partial}_{A'}s=0$. Now, Hartog theorem
implies that $s$ extends to a global section 
$\overline{s}\in L^p_2(X;E\times_{\Ad}\imag\klie)$.
By continuity $\overline{s}$ leaves $h(A,\Phi)$ fixed, and 
it is semisimple (for this we need to use that the eigenvalues of $s$
are constant). This contradicts the fact that $(A,\Phi)$ (and so $h(A,\Phi)$)
is semisimple. So $s$ does not fix $h(A,\Phi)$. Finally, to prove that
$\lambda(h(A,\Phi);s)>0$ we argue as in the proof of Theorem
\ref{corr}, using Lemma \ref{gradient}.
\end{pf}

\begin{lemma}
Fix a positive constant $C_B$. There exist positive constants
$C_1, C_2$ such that the following holds.
Let $g\in\GGG_G(X_0)=\Omega^0(X_0;E\times_{\Ad}G)$ be
such that $|g|_{{\log},C^0}\leq C_B |g|_{{\log},L^1}<\infty$. 
Then $|g|_{{\log},C^0}\leq C_1\Psi^c_{h(A,\Phi)}(g)+C_2$.
\label{fita3}
\end{lemma}
\begin{pf}
Since $h(A,\Phi)$ is analitically stable, 
given any $B>0$ there exist constants $C_1$ and $C_2$ such that for any
$s\in\MetB$ there is an inequality
\begin{equation}
\sup |s|\leq C_1\Psi^c_{h(A,\Phi)}(e^s)+C_2.
\label{proper}
\end{equation}
Thanks to the preceeding lemma, this inequality is valid not only for 
$s\in\MetB$, but also for any 
$$s\in\Met(C_B)=\{s\in L^p_2(X_0;E\times_{\Ad}\imag\klie)
|\ \|s\|_{C^0}\leq C_B\|s\|_{L^1}\},$$ 
as one can see tracing the proof of Lemma \ref{contradic}.
\end{pf}

\begin{lemma}
There a positive constant $C'$ such that 
for any $g_{\sigma,\chi}$ and $h$ and
for big enough (depending on $g_{\sigma,\chi}$ and $h$) $t>0$,
\begin{align}
|e^{t g_{\sigma,\chi}}h^{-1}|_{{\log},C^0}
&\leq C' |e^{t g_{\sigma,\chi}}h^{-1}|_{{\log},L^1} \notag \\
|e^{t g_{\sigma,\chi}}|_{{\log},C^0}
&\leq C' (|e^{t g_{\sigma,\chi}}h^{-1}|_{{\log},C^0}+1). \label{lanova}
\end{align}
\end{lemma}
\begin{pf}
This is a consequence of Lemma \ref{comparison} and the fact
that $X$ is compact (so $|h|$ and $|h^{-1}|$ are bounded functions
on $X$).
\end{pf}
By the properties of the integral of the moment map we have that
\begin{equation}
\Psi_{h(A,\Phi)}(e^{t g_{\sigma,\chi}} h^{-1})=
\Psi_{(A,\Phi)}(e^{t g_{\sigma,\chi}})-\Psi_{(A,\Phi)}(h).
\label{fita4}
\end{equation}
Now, putting $C_B=C'$ in Lemma \ref{fita3}, we conclude that 
\begin{align*}
t\sup|g_{\sigma,\chi}| &= |e^{t g_{\sigma,\chi}}|_{\log,C^0} 
\leq C'(|e^{t g_{\sigma,\chi}}|_{\log,C^0}+1) 
\qquad\text{ by (\ref{lanova})}\\
&\leq C'(C_1\Psi_{h(A,\Phi)}(e^{t g_{\sigma,\chi}} h^{-1})+C_2+1) 
\qquad\text{ by Lemma \ref{fita3}}\\
&\leq C_1'\Psi_{(A,\Phi)}(e^{t g_{\sigma,\chi}})+C_2'
\qquad\qquad\qquad\quad\text{ by (\ref{fita4}).}
\end{align*}
This implies, reasoning like in Theorem \ref{corr}, that 
$\lambda((A,\Phi);-\imag g_{\sigma,\chi})>0$.
By Lemma \ref{equivalencia}, this is equivalent to 
$T^c_{\Phi}(\sigma,\chi)>0$. Hence $(A,\Phi)$ is $c$-stable.

\begin{remark} When $F$ is a vector space the proof that existence
of solution implies stability is much easier if one uses the principle
that {\it curvature increases in subbundles} (see for example
\cite{Br2}). This is a consequence of the fact that the maximal weights
of a linear action of $K$ on a vector space are very simple
(see Lemma \ref{pesvect}) and, specially, that the maximal weight of
any element $s\in\Lie(\GGG_K)$ is constant along $\GGG_G$ orbits
in $\AAA^{1,1}\times\SSS$ (see section \ref{banfield}).
\end{remark}

\subsection{Uniqueness of solutions}
The proof is exactly as in the finite dimensional case:
it follows from the convexity of the integral of the moment map.

\subsection{Nonsimple pairs}
\label{nosimple}
The Hitchin--Kobayashi correspondence which
we have proved applies only to simple pairs $(A,\Phi)$. This
restriction is not always satisfied. As an example, suppose that there are
elements in the centre $Z=Z(\glie)$ of $G$ which leave $F$ fixed
(trivial example: $F$ equal to a point). Any element 
$z\in Z$ gives an element of the Lie algebra of the gauge group, which
we still denote by $z$. This element is semisimple and for any
$t$ the exponential $\exp(tz)$ fixes all connections in $\AAA$, and
by our assumption fixes also $\Phi$. In this situation, the pair
$(A,\Phi)$ is not simple.

When our group $G$ is $\GL(V)$, there is a standard way to solve
this problem. We assume that the whole center $Z$ leaves $\Phi$ fixed. 
We have to {\it split} the equation 
in the $Z$ part and in the $G/Z$ part. This is done as follows.
Define $\GGG_G^0$ to be the set of gauge transformation with
determinant pointwise equal to 1, and suppose that there 
are no semisimple elements in the Lie algebra of $\GGG_G^0$
which leave $(A,\Phi)$ fixed; under this assumption we can
find an element $g\in\GGG_G^0$ so that $g(A,\Phi)$ solves the 
trace-free part of the equation (observe that our proof applies 
to this situation); then Hodge theory gives a central element in 
$\GGG_G$ which, composed with $g$, solves the complete equation.

This idea applies for any reductive Lie group $G$. 
We just need to give a generalisation of the condition
of having determinant pointwise equal to 1 which we imposed
to the elements in $\GGG_G^0$. This is given by the following

\begin{lemma}
Let $G$ be a reductive Lie group. There exists $k\geq 1$
and a morphism $\phi:G\to(\CC^*)^k$ such that $\Ker\phi\cap Z$
is a discrete subgroup of $G$.
\end{lemma}
\begin{pf}
Take a faithful representation $\rho:G\to GL(W)$. Split $W$
in eigenspaces of the roots of $Z$ acting on $W$:
$W=W_1\oplus\dots\oplus W_k$, so that any central element
$z\in Z$ acts on any piece $W_j$ by homotecies.
Then $\rho(G)\subset GL(W_1)\times\dots\times GL(W_k)$, so
that for any $g\in G$ we have $\rho(g)=(g_1,\dots,g_k)$.
Let $\phi:G\to(\CC^*)^k$ be defined as
$\phi(g)=(\det g_1,\dots,\det g_k)$. Now suppose that there exists
$s\in Z(\glie)$ such that, for any $t$, 
$\phi(e^{ts})=(1,\dots,1)$. Since $e^{ts}$ acts
by homotecies on each piece, we must have 
$\rho(e^{ts})\in Z(SL(W_1))\times\dots\times Z(SL(W_k))
\simeq \ZZ/w_1\ZZ\times\dots\times\ZZ/w_k\ZZ$ for any $t$,
where $w_j=\dim W_j$. This implies that $\rho(e^{ts})=(1,\dots,1)$
and, since $\rho$ is faithful, $z=0$. This proves that 
$\Ker\phi\cap Z$ is discrete.
\end{pf}

Suppose now for simplicity that the whole center $Z(G)$ leaves
$\Phi$ fixed. 
We then define $\GGG_G^0$ to be the set of gauge transformations
which fibrewise belong to $\Ker\phi$, and proceed as in the case
$G=GL(V)$: we find $g\in\GGG^0_G$ such that the center free part
of the equation is solved and then use Hodge theory to solve the
complete equation.

\section{Yang-Mills-Higgs functional}
\label{YMHs}
In order to define the Yang-Mills-Higgs functional for pairs
in $\AAA^{1,1}\times\SSS$ it will be necessary to extend the
definition of covariant derivations on vector bundles to 
general fibre bundles. Recall that the subbundle $T\FFF_v$ of
vertical tangent vectors to $\FFF$ is by definition
$\Ker d\pi_F$, where $\pi_F:\FFF\to X$ is the projection. 
Using the Kaehler metric on $TF$ we get an induced metric on $T\FFF_v$
(recall that the action of $K$ respects the Kaehler structure and so in 
particular the Kaehler metric is kept fixed by $K$). 
In this section we will not use the fact that the complex structure
on $F$ is integrable, so that all the results remain valid when
$F$ is an almost-Kaehler manifold (in fact we could also consider
connections in $\AAA$).

\begin{definition}
Let $A\in\AAA^{1,1}$ be a connection on $E$. This connection 
induces a projection $\alpha:T\FFF\exh T\FFF_v$, since $\FFF$
is a fibre bundle associated to $E$. Take a section 
$\Phi\in\SSS=\Gamma(\FFF)$. We define the {\rm covariant
derivation of $A$ on $\Phi$} as
$$d_A\Phi=\alpha (d\Phi)\in\Omega^1(\Phi^*T\FFF_v).$$
\end{definition}

On the other hand,
since the complex structure $I_F$ on $F$ is left fixed by the action of
$K$, the bundle $\Phi^*T\FFF_v$ has an induced complex structure.
This justifies the following definition.

\begin{definition}
Let $(A,\Phi)\in\AAA^{1,1}\times\SSS$. We define the 
{\rm $\overline{\partial}$-operator of $A$ on $\Phi$} (resp. the 
{\rm $\partial$-operator of $A$ on $\Phi$} ) to be 
$\overline{\partial}_A\Phi=\pi^{0,1}d_A\Phi,$
(resp. $\partial_A\Phi=\pi^{1,0}d_A\Phi$), where
$\pi^{0,1}$ (resp. $\pi^{1,0}$) denotes the projection 
of $\Omega^1(\Phi^*T\FFF_v)$ to the second (resp. first) summand 
in the decomposition
$\Omega^1(\Phi^*T\FFF_v)=\Omega^{1,0}(\Phi^*T\FFF_v)
\oplus\Omega^{0,1}(\Phi^*T\FFF_v).$
\end{definition}

When $F$ is a vector space,
the operators $d_A$, $\partial_A$ and $\ov{\partial}_A$ coincide with 
the usual ones for vector bundles (in this case
there is a canonical identification $T\FFF_v\simeq\FFF$).

Recall that we have on $\klie$ a nondegenerate biinvariant positive 
definite pairing $\la,\ra$. This pairing gives a $K$-equivariant isomorphism 
$\klie\simeq\klie^*$ and an Euclidean metric on $\klie$ and $\klie^*$.

\begin{definition}
Fix a central element $c\in\klie$. 
The {\rm Yang-Mills-Higgs} 
functional $\YMH_c:\AAA^{1,1}\times\SSS\to\RR$ is defined as
$$\YMH_c(A,\Phi)=\|F_A\|_{L^2}^2+\|d_A\Phi\|_{L^2}^2
+\|c-\mu(\Phi)\|_{L^2}^2,$$
where $\Phi\in\SSS$ is a section and $A\in\AAA^{1,1}$ a connection on $E$.
\end{definition}

We will say that two sections $\Phi_0,\Phi_1\in\SSS$ are
{\rm homologous} if they induce the same map in cohomology, i.e.,
$\Phi_0^*=\Phi_1^*:H^*(\FFF)\to H^*(X).$

\begin{theorem}
\label{minimitza}
Fix a section $\Phi_0\in\SSS$.
The pairs $(A,\Phi)\in\AAA^{1,1}\times\SSS$ which minimize
the functional $\YMH_c$ among the pairs whose section
is homologous to $\Phi_0$ are those which satisfy the following pair
of equations
\begin{equation}
\left\{\begin{array}{l}
\overline{\partial}_A\Phi=0 \\
\Lambda F_A+\mu(\Phi)=c.\end{array}\right.
\label{vortex}
\end{equation}
\label{idminim}
\end{theorem}

The proof of this theorem will be given at the end of this section.

\subsection{}
The symplectic form $\omega_F$ gives an element
of $\Omega^0(\Lambda^2 (T\FFF_v)^*)$, since the action of $K$
keeps $\omega_F$ fixed. On the other hand, the connection $A$ on $E$ 
induces a projection $$\alpha:T\FFF\exh T\FFF_v$$ 
onto the subbundle of vertical tangent vectors. From this we obtain a map 
$\alpha^*:\Lambda^2 (T\FFF_v)^*\to \Lambda^2 T^*\FFF$, and we 
set $\tilde{\omega}_F^A=\alpha^*(\omega_F)\in\Omega^0(\Lambda^2 T^*\FFF)
=\Omega^2(\FFF)$. This 2-form is not in general closed. Consider the 2-form 
$\omega_F^A=\tilde{\omega}_F^A-\la \pi_F^*F_A,\mu\ra_{\klie}$. 

\begin{prop}
The 2-form $\omega_F^A\in\Omega^2(\FFF)$ is closed, and the cohomology class 
it represents is independent of the connection $A$.
\label{nonlinka2}
\end{prop}
\begin{pf}
The form $\omega_F^A$ coincides with the coupling form $\omega_{A,F}$ of
the symplectic fibration $\FFF\to X$ and the connection $A$ as defined
in \cite[Theorem 1.4.1]{GLeS}. This is proved in \cite[Example 2.3]{GLeS}.
In \cite[Theorem 1.4.1]{GLeS} it is proved that $\omega_F^A$
is closed and in \cite[Theorem 1.6.1]{GLeS} it is shown that the cohomology
of $\omega_F^A$ is independent of the connection $A$. 
\end{pf}

\begin{remark} One can prove that $\omega_F^A$ is the image by the 
generalised Chern-Weil homomorphism (see \cite[Chapter 7]{BeGeV})
of the equivariant de Rham form $\overline{\omega}_F=\omega_F-\mu$.
This gives another proof of Proposition \ref{nonlinka2} (see \cite{Mu}).
\end{remark}

In the sequel we will denote by $[\omega_{\FFF}]$ the cohomology class
represented by $\omega_F^A$. By a slight abuse of notation we will also
denote by $[\omega_{\FFF}]$ any de Rham form representing it.

\begin{prop}
For any section $\Phi\in\FFF$ and for any connection $A\in\AAA^{1,1}$, 
the following equality holds:
$$\int_X\la \Lambda F_A,\mu(\Phi)\ra_{\klie}=
\frac{1}{2}(\|\partial_A\Phi\|^2_{L^2}
-\|\overline{\partial}_A\Phi\|^2_{L^2})-
\int_X \Phi^*[\omega_{\FFF}]\wedge \omega^{[n-1]}.$$
\label{nonlinka}
\end{prop}

To prove Proposition \ref{nonlinka} we will use the following elementary 
lemma.

\begin{lemma}
Let $V$ and $W$ be two Euclidean vector spaces with scalar products
$\la,\ra_V$ and $\la,\ra_W$. Suppose that there are complex structures
$I_V\in\End(V)$, $I_W\in\End(W)$ and symplectic forms 
$\omega_V\in\Lambda^2V^*$, $\omega_W\in\Lambda^2W^*$ which satisfy
the following: $\la\cdot,\cdot\ra_V=\omega_V(\cdot,I_V\cdot)$
and $\la\cdot,\cdot\ra_W=\omega_W(\cdot,I_W\cdot)$.
Take a linear map $f:V\to W$ and let
$f^{1,0}$ (resp. $f^{0,1}$) be $(f+I_W\circ f\circ I_V)/2$
(resp. $(f-I_W\circ f\circ I_V)/2$). 
Let $2n=\dim_{\RR}V$. Then
$f^*\omega_W\wedge\omega_V^{[n-1]}=
\frac{1}{2}(|f^{1,0}|^2-|f^{0,1}|^2)\omega_V^{[n]},$
where, for any $g\in\Hom(V,W)$, $|g|^2=\Tr g^*g$.
\label{trucu}
\end{lemma}

\begin{pf} (Proposition \ref{nonlinka})
Using Lemma \ref{trucu} we have 
$$\int_X\Phi^*\tilde{\omega}_F^A\wedge\omega^{[n-1]}=
\frac{1}{2}(\|\partial_A\Phi\|^2_{L^2}
-\|\overline{\partial}_A\Phi\|^2_{L^2})$$
for any section $\Phi:X\to\FFF$. To apply the lemma we set, for
any $x\in X$, $V=T_xX$ and $W=T_{\Phi(x)}\FFF_v$ with the 
induced Kaehler structures, and $f=d_A\Phi(x)$. With these identifications
$f^{1,0}=\partial_A\Phi(x)$ and $f^{0,1}=\overline{\partial}_A\Phi(x)$. 
As a consequence, 
\begin{align}
&\frac{1}{2}(\|\partial_A\Phi\|^2_{L^2}
-\|\overline{\partial}_A\Phi\|^2_{L^2})-
\int_X \la \Lambda F_A,\mu(\Phi)\ra \notag \\
&=\int_X(\Phi^*\tilde{\omega}_F^A
-\Phi^*\la \pi_F^*\Lambda F_A,\mu(\Phi)\ra)\wedge\omega^{[n-1]}=
\int_X\Phi^*[\omega_{\FFF}]\wedge\omega^{[n-1]}.\notag 
\end{align}
This proves Proposition \ref{nonlinka}.
\end{pf}

\subsection{Proof of Theorem \ref{idminim}}

The following computation has its origins in an idea of 
Bogomolov in studying vortex equations on $\RR^2$. Here we
mimic \cite{Br1}, except that where he uses the Kaehler
identities we use Proposition \ref{nonlinka}.

\begin{lemma}
For any section $\Phi\in\SSS$ and any connection $A\in\AAA^{1,1}$ 
\begin{align*}
\YMH_c(A,\Phi)&=\|\Lambda F_A+\mu(\Phi)-c\|_{L^2}^2+
2\|\overline{\partial}_A\Phi\|_{L^2}^2
+2\int_X\la\Lambda F_A,c\ra \\
&+2\int_X \Phi^*[\omega_{\FFF}]\wedge \omega^{[n-1]}
-\int_X B_2(F_A,F_A)\wedge \omega^{[n-2]},
\end{align*}
where $B_2:\Omega^2(E\times_{\Ad}\klie)\otimes\Omega^2(E\times_{\Ad}\klie)
\to\Omega^4(X)$ denotes the combination of the wedge product with the
biinvariant pairing on $\klie$.
\end{lemma}
\begin{pf}
Throughout the proof $\|\cdot\|$ will denote $L^2$ norm. 
For any connection $A\in\AAA$ we have
$$|F_A|^2\omega^{[n]}=
|\Lambda F_A|^2\omega^{[n]}-B_2(F_A,F_A)\wedge\omega^{[n-2]}
+4|F_A^{0,2}|^2\omega^{[n]}$$
(see \cite[p. 209]{Br2}). 
We now develop using Proposition \ref{nonlinka} and
taking into account that $F_A^{0,2}=0$
\begin{align*}
&\|\Lambda F_A+\mu(\Phi)-c\|^2+2\|\overline{\partial}_A\Phi\|^2 \\
&= \|\Lambda F_A\|^2+\|\mu(\Phi)-c\|^2+
2\int_X\la\Lambda F_A,\mu(\Phi)\ra_{\klie}-2\int_X\la\Lambda F_A,c\ra 
+2\|\overline{\partial}_A\Phi\|^2\\
&= \|F_A\|^2+\|\mu(\Phi)-c\|^2+\|\partial_A\Phi\|^2+
\|\overline{\partial}_A\Phi\|^2
-2\int_X \Phi^*[\omega_{\FFF}]\wedge \omega^{[n-1]}
-2\int_X\la\Lambda F_A,c\ra \\
&+\int_X B_2(F_A,F_A)\wedge \omega^{[n-2]} \\
&= \|F_A\|^2+\|\mu(\Phi)-c\|^2+\|d_A\Phi\|^2
-2\int_X \Phi^*[\omega_{\FFF}]\wedge \omega^{[n-1]}
-2\int_X\la\Lambda F_A,c\ra \\
&+\int_X B_2(F_A,F_A)\wedge \omega^{[n-2]}.
\end{align*}
\end{pf}

Theorem \ref{minimitza} follows easily from the preceeding lemma. Indeed, 
$$2\int_X\la\Lambda F_A,c\ra+
2\int_X \Phi^*[\omega_{\FFF}]\wedge \omega^{[n-1]}
-\int_X B_2(F_A,F_A)\wedge \omega^{[n-2]}$$
is a topological quantity, that is, it only depends on the homology class
of $\Phi$. That this is true for the second summand is clear; as for
the first summand, by Chern-Weil theory one sees that it is equal
to a linear combination whose coefficients depend on $c$ of first
Chern classes of line bundles obtained from $E$ through representations
$K\to\UU$. Finally, the form $B(F_A,F_A)/8\pi^2$ represents the second
Chern character $ch_2\in H^4(X;\RR)$ of $V=E\times_{\rho_a}W_a$
(see \cite[p. 209]{Br2}); hence the third summand is also topological.

Finally, we obtain from \ref{minimitza} the following corollary
{\it {\`a} la} Bogomolov

\begin{corollary}
Suppose that a pair $(A,\Phi)$ is gauge equivalent to a pair satisfying
equations (\ref{vortex}). Then the following inequality holds
$$\int_X\la\Lambda F_A,c\ra+
\int_X \Phi^*[\omega_{\FFF}]\wedge \omega^{[n-1]}
-\frac{1}{2}\int_X B_2(F_A,F_A)\wedge\omega^{[n-2]}\geq 0.$$
\label{bogom}
\end{corollary}

\section{Example: the theorem of Banfield}
\label{banfield}
Suppose that $F$ is a Hermitian vector space and that $K$ acts on $F$
through a unitary representation $\rho:K\to U(F)$. D. Banfield \cite{Ba} has 
recently proved a general Hitchin--Kobayashi correspondence for this 
situation. The work of Banfield generalises existing results on vortex 
equations, Hitchin equations, and on other equations arising from particular 
choices of $K$ and $\rho$. In this section we will see how the result of 
Banfield can be deduced from Theorem \ref{main}. 

\subsection{The stability condition}
Let $h$ be the Hermitian metric on $F$. The
imaginary part of $h$ with reversed sign defines a symplectic form 
$\omega_F$ compatible with the complex structure and hence a Kaehler
structure. The action of $K$ on $F$ respects
the Kaehler structure and admits a moment map $\mu:F\to\klie^*$
$$\mu(x)=-\frac{\imag}{2}\rho^*(x\otimes x^*).$$ 
In other words, for any $s\in\klie$, 
$\la\mu(x),s\ra_{\klie}=-\frac{\imag}{2} h(x,\rho(s)x)$.
Let $x\in F$ and take an element $s\in\klie$.
Since $\rho(s)\in \ulie(F)$, the endomorphism $\rho(s)$ diagonalizes
in a basis $e_1,\dots,e_n$: $\imag\rho(s)e_k=\lambda_k e_k$,
where $\lambda_k$ is a real number for any $k$. 
Write $x=x_1e_1+\dots+x_ne_n$. 
\begin{lemma}
If $\lambda_k\leq 0$ for every $k$ such that $x_k\neq 0$,
then the maximal weight $\lambda(x;s)$ is equal to zero. Otherwise it is 
$\infty$.
\label{pesvect}
\end{lemma}

Let us assume that the representation $\rho$ is contained
in the representation $\rho_a$.
Let $E\to X$ be a $G$-principal bundle on a compact
Kaehler manifold $X$. Let $\FFF=E\times_{\rho} F$ be
the vector bundle associated to $E$ through the representation
$\rho$. Take a pair $(A,\Phi)\in\AAA^{1,1}\times\SSS$, and
fix a central element $c\in\klie$. Consider
on $E$ the holomorphic structure given by $\overline{\partial}_A$.
According to definition \ref{parella_estable},
$(A,\Phi)$ is $c$-stable if and only if for any parabolic subgroup
$P\subset G$, for any holomorphic reduction $\sigma\in\Gamma(X_0;E(G/P))$ 
defined on the complement of a complex codimension 2 submanifold $X_0$
of $X$ and for any antidominant character $\chi$ of $P$, the total degree
is positive:
$$T^c_{\Phi}(\sigma,\chi)>0.$$
The total degree is the sum of $\deg(\sigma,\chi)$ plus the
maximal weight of the action of $g_{\sigma,\chi}$ on 
$\Phi$ plus $\la\imag\chi,c\ra\Vol(X)$. The maximal weight is
\begin{equation}
\int_{x\in X}\lambda(\Phi(x);-\imag g_{\sigma,\chi}(x)).
\label{intBan}
\end{equation}
Define now $\FFF^-=\FFF^-(\sigma,\chi)\subset\FFF$ to be 
the subset given by the vectors in $\FFF$ on which $g_{\sigma,\chi}(x)$ 
acts negatively, that is, $v\in\FFF_x$ belongs to $\FFF^-$ if and only if
you can write $v=\sum v_n$ such that
$g_{\sigma,\chi}(x)(v_n)=\lambda_n v_n$ and $\lambda_n\leq 0$.
Since the eigenvalues of $g_{\sigma,\chi}$ are constant, $\FFF^-$
is a subbundle. And since the parabolic reduction is holomorphic,
so is $\FFF^-$.

If $\Phi\subset\FFF^-$, then the maximal weight at each fibre
is equal to zero by Lemma \ref{pesvect}, so the stability 
condition reduces to $\deg(\sigma,\chi)>0.$
On the other hand, if $\Phi(x)\notin\FFF^-_x$, then there is
an open neighbourhood $U$ of $x$ such that $\Phi(y)\notin\FFF^-_y$
for any $y\in U$. In this situation Lemma \ref{pesvect} tells us
that, for any $y\in U$, $\lambda(\Phi(y);-\imag g_{\sigma,\chi}(y))=\infty$.
Since this happens in an open set, the integral (\ref{intBan})
is infinite (since $X$ is compact, $\Phi$ is bounded and so 
$\lambda(\Phi(x);-\imag g_{\sigma,\chi}(x))$ is bounded below).
But the degree $\deg(\sigma,\chi)$ is always a finite number,
so the total degree will be positive (infinite, in fact) in this case. 
To sum up,

\begin{prop}
The pair $(A,\Phi)$ is stable if and only if for any
$P,\sigma,\chi$ as above, if $\Phi$ is contained in 
$\FFF^-(\sigma,\chi)$, then $\deg(\sigma,\chi)+\la\imag\chi,c\ra\Vol(X)>0.$
\end{prop}

This is precisely Banfield stability condition.

\subsection{Simple pairs}
To give a characterisation of simple pairs we use the following
definition due to Banfield \cite{Ba}:

\begin{definition}
Suppose that the vector bundle $\FFF$ decomposes into a nontrivial
direct sum $\bigoplus_k \FFF_k$ of holomorphic vector bundles
and that there is a reduction of the structure group of $E$
to $G'\subset G$, compatible with the splitting. 
Suppose further that a central element of the Lie algebra 
$\glie'$ of $G'$ annihilates the section $\Phi$ but acts nontrivially 
on $\FFF$. Then we say $(A,\Phi)$ is a {\rm decomposable pair}. 
If no such splitting exists, the we say that $(A,\Phi)$ is an 
{\rm indecomposable pair}.
\end{definition}

\begin{lemma}
The pair $(A,\Phi)$ is simple if and only if 
it is indecomposable.
\label{bansim}
\end{lemma}
\begin{pf}
Suppose that 
$0\neq s\in\Omega^0(E\times_{\Ad}\glie)$ is semisimple and stabilises 
$(A,\Phi)$. In particular $\fX^{\AAA}_s(A)=0$, and this implies that 
$\overline{\partial}_A(s)=0$. So the eigenvalues of $\rho(s)$ 
are constant, and since $s$ is semisimple $\rho(s)$ diagonalises. 
Let the different eigenvalues of $\rho(s)$ be $\lambda_1<\dots<\lambda_r$, 
and consider the decomposition 
$\FFF=\FFF(\lambda_1)\oplus\dots\oplus\FFF(\lambda_r)$
in eigenbundles, which are holomorphic,
and every $\FFF_k=\FFF(\lambda_k)$ having as structure group
a subgroup $G_k\subset G$. 
Since $s$ leaves $\Phi$ fixed
$\Phi$ must belong to $\FFF(0)$. On the other hand, $0$ in obviously
not the unique eigenvalue of $\rho(s)$, so the decomposition
$$\FFF=\FFF_1\oplus\dots\oplus\FFF_r$$
is not trivial. Finally, the section $s$ provides the central element
killing $\Phi$.

The proof of the converse is similar.
\end{pf}

\subsection{The equations}
Our equation (\ref{hk}) in the case of linear representations
is the same one given by Banfield (note that Banfield also
considers the holomorphicity condition $\ov{\partial}_A\Phi=0$). 

\section{Example: filtrations of vector bundles}
\label{filtracions}

In this section we study Theorem \ref{main} in the particular case
in which $F$ is a Grassmannian or, more generaly, a flag manifold.
We assume, for simplicity, that $X$ is a Riemann
surface. For the higher dimensional case everything that follows
remains valid if we consider reflexive subsheaves and not only
subbundles in the definition of stability (this reflects the need of
considering reductions of the structure group defined on the
complement of a complex codimension 2 submanifold of $X$ in
the general definition of stability).

\subsection{Projective manifolds with actions of Lie groups} 
Let $F\subset\PP(\CC^n)$ be any smooth complex subvariety.
Let us take on $\CC^n$ the canonical Hermitian metric.
This allows to define on $\PP(\CC^n)$ the Fubini-Study Kaehler 
structure. We consider on $F$ the induced structure.
Suppose that a compact Lie group $K$ acts on $\PP(\CC^n)$ 
through a representation $\rho:K\to \U(n;\CC)$
leaving $F$ fixed. Since $\rho(K)\subset \U(n;\CC)$,
the action of $K$ on $\PP(\CC^n)$ (and hence on $F$) respects the 
Kaehler structure. A moment map $\mu_F:F\to\klie^*$ for this action is
\begin{equation}
\mu_F(x)=-\frac{\imag}{2}\rho^*\left(\frac{\hat{x}\otimes \hat{x}^*}
{\|\hat{x}\|^2}\right),
\label{momproj}
\end{equation}
where $\hat{x}\in \CC\setminus\{0\}$ denotes any lift of $x\in F$. 
Take a point $x\in F$ and consider an element $s\in\klie$. We can take a 
basis $e_1,\dots,e_n$ of $\CC^n$ in which the action of $s$ diagonalizes:
for any $k$, $\imag\rho(s)e_k=\lambda_k e_k$, where $\lambda_k$ is 
a real number. Fix a lifting $\hat{x}\in \CC\setminus\{0\}$ of $x\in F$ 
and write $\hat{x}=x_1e_1+\dots+x_ne_n$. 
\begin{lemma}
The maximal weight of $s$ acting on $x$ is 
$\lambda(x;s)=\max \{\lambda_k|x_k\neq 0\}.$
\label{pesmaxproj}
\end{lemma}

The manifold $F$ will be in this section either a Grassmannian
or a flag manifold. The Lie group $K$ will be $\U(R;\CC)$,
where $R\geq 1$ is an arbitrary integer, and
we will take the standard representation in $\CC^R$ as our
representation. We will assume for simplicity that
$\Vol(X)=1$.

\subsection{Subbundles}
Let $E\to X$ be a principal $\U(R;\CC)$ bundle on $X$. Consider
the standard representation on $\CC^R$. The associated bundle is 
a vector bundle $V\to X$ of rank $R$. Using Theorem \ref{main}, we 
will find a Hitchin--Kobayashi correspondence for subbundles 
$V_0$ of $V$ of fixed rank $0<k<R$. This correspondence has already
been proved in \cite{BrGP1} and in \cite{DaUW}.

Using an idea of \cite{DaUW} we identify the inclusion 
$V_0\hookrightarrow V$ with a section $\Phi$ of the bundle with fibres the 
Grassmannian of $k$-subvectorspaces $\Gr_k(\CC^R)$ associated to
$E$ by the usual action of $\GL(R;\CC)$ on $\Gr_k(\CC^R)$:
$$\FFF=E\times_{\GL(R;\CC)} \Gr_k(\CC^R).$$
The Pl{\"u}cker embedding maps $\Gr_k(\CC^R)$ in a $\GL(R;\CC)$-equivariant 
way into $\PP(\Lambda^k\CC^R)$, and the action 
of $\GL(R;\CC)$ in $\PP(\Lambda^k\CC^R)$ lifts to the obvious
action in $\Lambda^k\CC^R$. So we are in the situation described
at the beginning of this section.
Observe that the centre of $\GL(R;\CC)$ acts trivially on the Grassmannian. 
In consequence, the comments in subsection \ref{nosimple} are relevant in
this situation.

If $\omega$ is the symplectic form in $\Gr_k(\CC^R)$ 
inherited by the Fubini-Study symplectic form on $\PP(\Lambda^k\CC^R)$, 
then $\tau\omega$ also gives $\Gr_k(\CC^R)$ a Kaehler structure when 
$\tau>0$ and everything gets multiplied by $\tau$: the moment map, the
maximal weights and the integral of the moment map.
We fix from now on a constant $\tau>0$ and
we work with the symplectic form $\tau\omega$. The constant
$\tau$ can be identified with the parameter appearing in the
notion of stability and in the equations in \cite{BrGP1,DaUW}.

\subsection{The moment map for the Grassmannian with the action of $\U(n)$}
The action of $\U(n;\CC)$ on $\Gr_k(\CC^R)$ is symplectic. 
Making use of formula (\ref{momproj}) one easily verifies
that if $\pi\in\Gr_k(\CC^R)$, then the moment map of the action 
of $\U(n;\CC)$ at the point $\pi$ is the element in $\ulie(n;\CC)^*$
which sends $\xi\in\ulie(n;\CC)$ to
$\mu(\pi)(\xi)=-\imag\tau\Tr(\pi\circ\xi)$,
where $\pi$ denotes the orthogonal projection onto $\pi$ (see
\cite{DaUW}, p. 485).

\subsection{Maximal weights of $\U(n)$ acting on the Grassmannian}
Consider the standard action of $\U(n)$ on $\PP(\Lambda^k\CC^R)$.
Take an element $s\in\ulie(n)$.
We now give the maximal weight $\lambda(v;s)$ in the
case when $v=v_1\wedge\dots\wedge v_k\neq 0$, for $v_j\in \CC^R$.
This case is enough for our purposes, since the image of the Grassmanian
$\Gr_k(\CC^R)$ given by the Pl{\"u}cker embedding into $\Lambda^k\CC^R$ 
is precisely the set of points of that form. 

Let $\pi$ be the $k$-subspace of $\CC^R$ spanned by $\{v_j\}$.
Let $\lambda_1<\dots<\lambda_r$ be the eigenvalues of
$\imag s$ acting on $\Lambda^k\CC^R$, 
and for any $1\leq j\leq r$ write 
$E_j=\bigoplus_{i\leq j}\Ker(\imag s-\lambda_k\Id)$.
Set $\alpha_j=\lambda_j-\lambda_{j+1}$. Then
\begin{equation}
\lambda(v;s)=\tau\left(
\dim(\pi)\lambda_r+\sum_{j=1}^{r-1}\dim(\pi\cap E_j)\alpha_j\right).
\label{grass}
\end{equation}
The proof of this formula is an easy exercise which follows
from Lemma \ref{pesmaxproj}.

\subsection{Simple extensions}
Reasoning similarly as in Lemma \ref{bansim} one can prove this
\begin{lemma}
The pair $(A,\Phi)$ is not simple if and only if one can find
a holomorphic (with respect to $\overline{\partial}_A$)
splitting $V=V'\oplus V''$ such that the subbundle $V_0$
given by the section $\Phi$ is contained in $V'$.
\end{lemma}

\subsection{The stability condition}

Let $c\in\RR$ be a real number. Fix a pair $(A,\Phi)$,
which gives a holomorphic structure on $V$ and an inclusion of bundles
$V_0\subset V.$ 
In this section we will study the $-\imag c\Id$-stability condition
for the pair in terms of $V_0\subset V$.

A (holomorphic) parabolic reduction $\sigma$ of the structure 
group of $E$ is the same as giving a (holomorphic) filtration 
$0\subset V^1\subset\dots\subset V^{r-1}\subset V^r=V,$
and an antidominant character $\chi$ for this reduction is of the form
$\chi=z\Id+\sum_{j=1}^{r-1}m_j\lambda_{R^j},$
where $R^j=\rk(V^j)$, $\lambda_{R^j}=\pi_j-\frac{R^j}{R}\Id$,
$\pi_j$ is the projection onto $\CC^{R^j}$),
$z$ is any real number, and the $m_j$ are real negative numbers.
Taking into account that the representation is just
the standard representation of $\GL(n;\CC)$ in $\CC^R$ we deduce that the                   
degree of the pair $(\sigma,\chi)$ is
$$\deg(\sigma,\chi)=z\deg(V)+\sum_{j=1}^{r-1}m_j
\left(\deg(V^j)-\frac{R^j}{R}\deg(V)\right).$$

To calculate the maximal weight of the action of $\chi$ on the section
$\Phi$ we use formula (\ref{grass}). The parameters that appear there
are related to ours as follows: $\alpha_j=m_j$ for any $1\leq j\leq r-1$
and $\lambda_r=z-\sum_{j=1}^{r-1}m_j\frac{R^j}{R}$.
We get, after integration
(recall that the volume of $X$ has been normalized to 1):
\begin{equation}
\int_{x\in X}
\mu(\Phi(x);-g_{\sigma,\chi}(x))=
\rk(V_0)\left(z-\sum_{j=1}^{r-1}m_j\frac{R^j}{R}\right)
+\sum_{j=1}^{r-1}m_j\rk(V_0\cap V^j).
\end{equation}
Hence, the stability notion is as follows: for any filtration
$0\subset V^1\subset\dots\subset V^{r-1}\subset V^r=V$ and any set
of negative weights $\alpha_1,\dots,\alpha_{r-1}$ we must have
\begin{align}
0 &< z\deg(V)+\sum_{j=1}^{r-1}m_j
\left(\deg(V^j)-\frac{R^j}{R}\deg(V)\right) \notag \\
  &+ \tau\left(\rk(V_0)
\left(z-\sum_{j=1}^{r-1}m_j\frac{R^j}{R}\right)
+\sum_{j=1}^{r-1}m_j\rk(V_0\cap V^j)\right)-zc\ R.
\label{st1}
\end{align}
(Observe that thanks to our assumption that $\Vol(X)=1$,
$\la\imag\chi,c\ra\Vol(X)=-zcR$.)
If this is to be satisfied by all possible choices of $z$, then
$$c=\frac{\deg(V)+\tau\rk(V_0)}{R}.$$
So, given the symplectic form $\tau\omega$, there is a unique central 
element $c\in\ulie(n;\CC)$ such that the pair can be
$c$-stable.
Putting the value of the central element inside (\ref{st1}) we get
\begin{align}
0 &< \sum_{j=1}^{r-1} m_j\left(
\deg(V^j)-\frac{R^j}{R}\deg{V}-\tau\rk(V_0)
\frac{R^j}{R}+\tau\rk(V_0\cap V^j)\right) \notag \\
&= \sum_{j=1}^{r-1} m_jR^j\left(
\frac{\deg(V^j)+\tau\rk(V_0\cap V^j)}{R^j}
-\frac{\deg(V)+\tau\rk(V_0)}{R}\right),\notag 
\end{align}
and using the fact that 
the numbers $m_j$ are arbitrary negative numbers, we see that
a necessary and sufficient condition for $(E,\Phi)$ to be stable
is that for any nonzero proper subbundle (in fact, reflexive subsheaf)
$V^1\subset V$
$$\frac{\deg(V^1)+\tau\rk(V_0\cap V^1)}{\rk(V^1)}
<\frac{\deg(V)+\tau\rk(V_0)}{R},$$
and this is the same condition that appears in \cite{DaUW, BrGP1}. 

In what concerns the equations, they are exactly those in \cite{DaUW}. 
Instead of writing them in terms of a gauge
transformation, we will put as the variable a metric $h$ in the 
bundle $V$. This is equivalent to our setting, since the relevant
space in our case is the gauge group of complex transformations modulo
unitary gauge transformations, and this coset space can be identified
with the space of metrics. Taking into account the precise form of the
moment map for the action of $\GL(n;\CC)$ in $\Gr_k(\CC^R)$
we can write the equations as $\Lambda F_A-\imag\tau\pi^h_{V_0}=-\imag c\Id,$
where $\pi^h_{V_0}$ is the $h$-orthogonal projection onto $V_0$. 
The equations considered in \cite{BrGP1} are written in a different way, 
but in \cite{DaUW} it is proved that they are equivalent to ours. 

\subsection{Filtrations}
Here we generalise the preceeding results to
the case of filtrations (see \cite{AlGP}).
Our trick is to identify a filtration
$0\subset V_1\subset\dots\subset V_s\subset V$
with a section $\Phi$ of the associated bundle with fibre the flag 
manifold $F_{i_1,\dots,i_s}$, where $i_k=\rk(V_k)$.
This manifold is embedded in a product of Grassmannians.
The Kaehler structure in the flag manifold is not unique.
We can in fact take as symplectic form any weighted sum 
of the pullbacks of the symplectic forms in the Grassmannians,
provided the weights are positive. So the Kaehler structure
depends on a $s$-uple of positive parameters $\tau=(\tau_1,\dots,\tau_s)$.
We can now work out the stability notion analogously to the case of
extensions, and obtain that (here we write 
$0\subset V_1\subset\dots\subset V_s\subset V$ for the filtration
represented by the section $\Phi$)
\begin{itemize}
\item the equation is 
$\Lambda F_A-\imag\sum\tau_k\pi^h_{V^k}=-\imag c\Id$, 
where $\pi^h_{V^k}$ is the $h$-orthogonal projection onto
$V^k$ and where $c$ is a real constant;
\item the pair $(A,\Phi)$ is simple unless there exists a holomorphic
(with respect to $\ov{\partial}_A$) splitting 
$V=V'\oplus V''$ such that $V_k\subset V'$ for any $k\leq s$;
\item the only value of $c$ for which we can expect our filtration
to be $c$-stable is 
$$c=\frac{\deg(V)+\sum\tau_k\rk(V_k)}{R};$$
\item the stability notion is as follows: for any nonzero proper
reflexive subsheaf $V^1\subset V$,
$$\frac{\deg(V^1)+\sum\tau_k\rk(V_k\cap V^1)}{\rk(V^1)}<
\frac{\deg(V)+\sum\tau_k\rk(V_k)}{R}.$$
\end{itemize}

\subsection{Bogomolov inequality}
In this subsection we state the Bogomolov inequality given in 
Corollary \ref{bogom} for the case of filtrations. 
For that we need to compute the cohomology class
$\Phi^*\phi_A(\ov{\omega}_F)$.

We begin with some general observations.
When the cohomology class represented by the symplectic form 
$\omega_F$ of $F$ belongs to $H^2(F;\imag 2\pi\ZZ)$, there
exists a line bundle $L\to F$ with a connection $\nabla$ whose curvature
coincides with $-\imag\omega_F$. Assume that the action of
$K$ on $F$ lifts to a linear action on $L$. Then $\nabla$ can
be assumed to be $K$-equivariant (by just averaging if it is not).
Using the action of $K$ on $L$ we can define a line bundle $\LLL\to\FFF$
as $\LLL=E\times_K L$. Denote $\plx:\LLL\to X$ and
$\plf:\LLL\to \FFF$ the projections. Let $A$ be a connection on $E$. 
The connection $A$ induces a connection on the associated bundle
$\LLL$, which may be seen as a projection 
$\alpha:T\LLL\to\Ker d\plx$. Since $\nabla$ is $K$-equivariant,
we may extend it fiberwise to obtain a projection
$\beta:\Ker d\plx \to\Ker d\plf$. The composition 
$\gamma=\beta\circ\alpha:T\LLL\to \Ker d\plf$ defines a connection
$\nabla^A$ on $\LLL\to\FFF$. It is an exercise to verify that 
$\omega_F^A=\imag F_{\nabla^A}$,
where $F_{\nabla^A}$ is the curvature of $\nabla^A$.

If $F=\Gr_k(\CC^R)$ is a Grassmannian everything in the preceeding
paragraph works. In particular, the line bundle $L\to F$ can
be identified with the dual of the determinant bundle, that is,
with the line bundle whose fiber on $V\in\Gr_k(\CC^R)$ 
is $\Lambda^k V^*$. More generaly, if $F=F_{i_1,\dots,i_s}$ and $F$ has the 
Kaehler structure induced by the parameters $\tau=(\tau_1,\dots,\tau_s)$, 
then for any $(A,\Phi)\in\AAA^{1,1}\times\SSS$ we have
$$\int_X\Phi^*[\omega_{\FFF}]\wedge\omega^{[n-1]}
=-\sum_{k=1}^s \tau_k\deg(V_k),$$
where $V_1\subset\dots\subset V_s\subset V$ is the filtration
represented by the section $\Phi$.

So Corollary \ref{bogom} takes the following form in this case:

\begin{corollary}
Let $A$ be a connection on $E$, and consider a filtration
$0\subset V_1\subset\dots\subset V_s\subset V$ which is holomorphic
with respect to $\ov{\partial}_A$. Let us write $\Phi$
for the section of $\FFF$ which represents this filtration.
If the pair $(A,\Phi)$ is $\GGG_G$ equivalent to a solution of 
$$\Lambda F_A-\imag\sum\tau_k\pi^h_{V^k}=-\imag c\Id,$$
then the following holds
$$\deg(V)\left(\frac{\deg(V)+\sum\tau_k\rk(V_k)}{R}\right)
-\sum_{k=1}^s\tau_k\deg(V_k)-4\pi^2 \la ch_2(V)
\cup\omega^{[n-2]},[X]\ra\geq 0.$$
\end{corollary}

\end{document}